\documentclass[11pt]{amsart}
\usepackage{amsmath,amssymb,amsthm}
\usepackage{hyperref}
\usepackage{mathrsfs}
\usepackage{graphicx}
\usepackage{tikz}
\usepackage{enumerate}
\usepackage{float}

\theoremstyle{plain}
\newtheorem{theorem}{Theorem}[section]
\newtheorem{lemma}[theorem]{Lemma}

\newtheorem{proposition}[theorem]{Proposition}

\theoremstyle{definition}

\newtheorem{remark}[theorem]{Remark}
\newtheorem{example}[theorem]{Example}

\newtheorem*{namedthm*}{\thistheoremname}
\newcommand{\thistheoremname}{} 

\newcommand{\da}{\mathord{\downarrow}}
\newcommand{\rom}[1]{\rm{\uppercase\expandafter{\romannumeral #1}}}

\newcommand{\II}{\begin{enumerate}}
\newcommand{\III}{\end{enumerate}}

\makeatletter
\def\ps@pprintTitle{%
  \let\@oddhead\@empty
  \let\@evenhead\@empty
  \def\@oddfoot{\reset@font\hfil\thepage\hfil}
  \let\@evenfoot\@oddfoot
}
\makeatother

\begin{document}

\title{Products of two sober dcpo's need not be  sober} \thanks{This research is supported by NSFC (No.\,12231007 and No.\,12371457).}

\author{Hualin Miao, Xiaoyong Xi, Xiaodong Jia, Qingguo Li, Dongsheng Zhao}
\address{H.\,Miao, School of Mathematics, Hunan University, Changsha, Hunan, 410082, China. Email: {\rm miaohualinmiao@163.com} }
\address{X.\,Xi, School of Mathematics and Statistics, Yancheng Teachers University, Yancheng 224002, China. Email: {\rm xixy@yctu.edu.cn}}
\address{X.\,Jia (corresponding author),  School of Mathematics, Hunan University, Changsha, Hunan, 410082, China. Email:  {\rm jiaxiaodong@hnu.edu.cn} }
\address{Q.\,Li, School of Mathematics, Hunan University, Changsha, Hunan, 410082, China.  Email: {\rm liqingguoli@aliyun.com}}
\address{D.\,Zhao,  Mathematics and Mathematics Education, National Institute of Education, Nanyang Technological University, 1 Nanyang Walk, 637616, Singapore. Email:  {\rm dongsheng.zhao@nie.edu.sg} }

\begin{abstract}
We  construct  two dcpo's whose Scott spaces are sober, but  the Scott space of their order product is not sober. This answers an open problem on the sobriety of Scott spaces. \textcolor{black}{In addition}, we show that if $M$ and $N$ are special  type of  sober complete lattices, then the Scott space of their order product $M\times N$ is sober.
\end{abstract}

\maketitle

\section{Introduction}

Sobriety of  topological spaces is one of the most extensively studied properties either in classical topology or in non-Hausdorff topology. In the former,  sobriety has already been used in characterizing the spectral spaces of commutative rings \cite{hochster1969prime}, while in non-Hausdorff topology,  sobriety plays a  key role in linking topological theory to locale theory through the celebrated Stone Duality \cite{goubault13a}. In domain theory, which was introduced by Dana Scott to give denotational semantics to programming languages, the study of the sobriety of Scott spaces has also got a relatively long history. For instance, Jimmie  Lawson and Karl Hofmann independently proved that the Scott topology of every domain (continuous directed complete poset) is sober \cite{gierz03}, and Samson Abramsky used the Stone duality to link Gordon Plotkin's domain-theoretic semantics (SFP-domains) to the logical descriptions of programs, along the lines of his celebrated Domain Theory in Logical Form~\cite{abramsky87a}.

Sobriety of topological spaces is relatively transparent, and the category of all sober spaces and continuous maps is a reflective subcategory of all ($T_0$) topological spaces, see e.\,g., \cite{gierz03}. However, sobriety on semantic domains or in general on directed complete posets (dcpo's for short) is quite opaque.  A dcpo that has a sober Scott topology is \textcolor{black}{often referred to as a sober dcpo}. In fact, it was unknown for a long time whether every dcpo is sober, and it was Peter Johnstone who constructed the first dcpo with a non-sober Scott topology \cite{johnstone81}. Johnstone's example has then inspired many further work. Based on Johnstone's work, John Isbell gave a complete lattice that is not a sober space in the Scott topology \cite{isbell82}, and remarkably recent years have witnessed \textcolor{black}{a lot of research progress} in non-Hausdorff topology that is stimulated by sobriety \cite{Miaoxl2023, xuxizhao2021}. 

As mentioned above, sober spaces reflect among all $T_0$ topological spaces, and it follows easily from category theory that sober spaces are closed under taking topological products. However, a similar result is not known for sober dcpo's. The obstacle is hidden in the both angelic and demonic fact that the Scott \textcolor{black}{topology of the product} of two dcpo's is in general different from the product of the Scott topology on the two dcpo's. 
 Indeed, the third author asked questions about forms of irreducible closed sets in  products of sober dcpo's \cite{jia2020order}, and in their study on the dcpo's determined by their lattices of Scott topology in~\cite{zhao2018uniqueness}, the authors wish to know whether the product of two sober dcpo's is sober. Although this problem was first explicitly stated  in~\cite{xu2021some}, it has been open in the community of domain theory for a long time. 

The main objective of this paper is to give a negative answer to the above problem. We shall construct two sober dcpo's $P_1$ and $P_2$ such that $P_1\times P_2$ is not sober in the Scott topology. Hence, the category of all sober dcpo's is not reflective in the category of all dcpo's with all Scott continuous mappings as morphisms. 
We will also explore sober complete lattices of a  special  type, and then prove that their finite products are indeed sober in the Scott topology. 


\section{Preliminaries}
In this section, we recall some basic concepts and notations that will be used in this paper.

Let $P$ be a partially ordered set (\emph{poset}, for short), $D\subseteq P$ is \emph{directed} (resp., \emph{filtered}) if $D$ is nonempty and for any finite subset $F\subseteq D$, there is $d\in D$ such that $d$ is an upper bound (resp., a lower bound) of $F$. A poset $P$ is called \emph{directed complete} (\emph{dcpo}, for short) if every directed subset $D$ of $P$ has a least upper bound (supremum), which we denote by $\sup D$, or $\bigvee D$. \textcolor{black}{It is easy to see that a poset $P$ is a dcpo if and only if every directed subset of $P$, which does not have a maximum element, has a supremum in $P$. We call directed subsets that do not have maximum elements \emph{non-trivial directed subsets}}. 
For any subset $A\subseteq P$, let $\uparrow$$A$ = $\{x\in P: x\geq$ a for some $a\in A\}$ and $\downarrow$$A$ = $\{x\in P: x\leq a$ for some $a\in A\}$. Specifically, we write $\uparrow$$x$ = $\uparrow$$\{x\}$ and $\downarrow$$x$ = $\downarrow$$\{x\}$. We will call $A\subseteq P$ an \emph{upper set} (resp., a \emph{lower set}) if $A$ = $\uparrow$$A$ (resp., $A$ = $\downarrow$$A$). A complete lattice is a poset in which every subset has a supremum and an infimum.

A subset $U$ of $P$ is \emph{Scott open} if $U$ = $\uparrow$$U$ and for any directed subset $D$ for which sup$D$ exists, sup$D$ $\in U$ implies $D\cap U \neq \emptyset$. Accordingly, $A\subseteq P$ is \emph{Scott closed} if $A$ = $\downarrow$$A$ and for any directed subset $D$ of $P$ with sup$D$ existing, $D\subseteq A$ implies sup$D$ $\in A$. The set of all Scott open sets of $P$ forms the \emph{Scott topology} on $P$, which is denoted by $\sigma (P)$, and the set of all Scott closed sets of $P$ is denoted by $\Gamma(P)$. Furthermore, for a subset $A$ of $P$, we will use $\overline{A}$ or $cl(A)$ to denote the closure of $A$ with respect to the Scott topology on $P$. The space $(P, \sigma (P))$ is also denoted by $\Sigma P$, some authors \textcolor{black}{call} such a space, a poset endowed with the Scott topology, a \emph{Scott space}.

For a $T_{0}$ space $X$, we denote all open sets of $X$ by $\mathcal{O}(X)$, the partial order $\leq$$_{X}$, defined by $x\leq_{X} y$ if and only if $x$ is in the closure of $y$, is called the \emph{specialization order} on $X$. Naturally, the closure of a single point $x$ is $\downarrow$$x$, the order considered here is, of course, the specialization order. The specialization order on a Scott space $\Sigma P$ coincides with the original order on $P$. A $T_0$ space $X$ is \emph{sober} if every nonempty irreducible closed subset $C$ of $X$ is the closure of some unique singleton set $\{c\}$, where $C\subseteq X$ is called \emph{irreducible} if $C\subseteq A\cup B$ for closed subsets $A$ and $B$ implies that $C\subseteq A$ or $C\subseteq B$. We denote the set of all closed irreducible subsets of $X$ by~$\mathrm{IRR}(X)$. 
\section{One sufficient condition}

In this section, we prove a positive result for the sobriety of products of dcpo's.

\begin{lemma}\label{product}
	Let $L, P$ be two countable dcpo's and $L^*=\sigma (L)$ and $P^*=\sigma (P)$ be the complete lattices of all Scott open sets of $L$ and $P$, respectively. Then the product topology of $\Sigma L^*$ and $\Sigma P^*$ coincides with the Scott topology of the product poset $L^*\times P^*$. That is, $\Sigma (L^*\times P^*)= \Sigma L^*\times \Sigma P^*.$
	\begin{proof}
		One sees immediately that the Scott topology of the product poset $L^*\times P^*$ \textcolor{black}{is at least as fine as} the corresponding product topology. For the converse, assume that there is a $\mathcal U\in \sigma(L^*\times P^*)$ such that $\mathcal U$ is not open in $\Sigma L^*\times \Sigma P^*$. \textcolor{black}{
       Recall sets of the form $\mathcal U'\times \mathcal V'$, with $\mathcal U'\in \sigma(L^*)$ and $\mathcal V'\in \sigma(P^*)$, form a basis for the product topology.
       So there exists $(U,V)\in \mathcal{U}$ such that 
for any $\mathcal{U'}\times \mathcal{V'}$ that contains $(U, V)$, we have $\mathcal{U'}\times \mathcal{V'}\not\subseteq \mathcal{U}$.
 From the countability of $L$ and $P$, \textcolor{black}{we can in addition assume} that $U=\{x_n:n\in \mathbb{N}\}$ and $V=\{y_n:n\in \mathbb{N}\}$. Now for each $n\in \mathbb{N}$, let 
		\begin{center}
			$\mathcal{U}_n=\{W\in \sigma(L):\{x_i:i\leq n\}\subseteq W\}$ and $\mathcal{V}_n=\{W\in \sigma(P):\{y_i:i\leq n\}\subseteq W\}$.
		\end{center}
		Obviously, $(U,V)\in \mathcal{U}_n\times \mathcal{V}_n$ and $\mathcal{U}_n\times \mathcal{V}_n$ is open in the product topology of $L^*\times P^*$  for any $n\in \mathbb{N}$. Hence, it follows that $\mathcal{U}_n\times \mathcal{V}_n\not\subseteq \mathcal{U}$ for all $n\in \mathbb{N}$.
		Pick $(W_n,W_n')\in (\mathcal{U}_n\times \mathcal{V}_n)\backslash \mathcal{U}$ for each $n\in \mathbb{N}$.}
		
		For each $n\in \mathbb{N}$, write
		\begin{center}
			$U_n=(\bigcap_{k\geq n}W_k)\cap U$, $V_n=(\bigcap_{k\geq n}W_k')\cap V$.
		\end{center}
		
		$\mathbf{Claim~1}$: For any $n\in \mathbb{N}$, $U_n$ and $V_n$ are both Scott open in $L$, $P$, respectively.
		
		It is easy to see that $U_n$ is an upper set of $L$ for any $n\in \mathbb{N}$. Let $D$ be a directed subset of $L$ with $\sup D$ existing and $\sup D\in U_n$. It suffices to check that $D\cap U_n\neq \emptyset$. We will only discuss the non-trivial case: $\sup D\notin D$. Note that $\sup D\in U$. Then we have some $d_0\in D$ such that $d_0\in U$ from the Scott openness of $U$. It follows that $d_0=x_{n_0}$ for some fixed $n_0\in \mathbb{N}$.
		
		If $n\geq n_0$, then the fact that $W_k\in \mathcal{U}_k$ implies that $x_{n_0}\in W_k$ for any $k\geq n$. This means that $x_{n_0}\in D\cap U_n$, that is, $D\cap U_n\neq \emptyset$. If $n\leq n_0-1$,  we notice that $\sup D\in \bigcap_{n\leq k\leq n_0-1}W_k$. Then there exists $d_1\in D\cap \bigcap_{n\leq k\leq n_0-1}W_k$ by the Scott openness of $\bigcap_{n\leq k\leq n_0-1}W_k$. We can find $d\in D$ such that $d_0,d_1\leq d$ because $D$ is directed. This yields that $d\in D\cap U_n$ since $d_0\in\bigcap_{k\geq n_0}W_{k}\cap U$ and $d_1\in \bigcap_{n\leq k\leq n_0-1}W_k$. Hence, $U_n$ is Scott open. By using similar deduction to the above proof, we know that $V_n$ is Scott open in $P$ for any $n\in \mathbb{N}$.
		
		$\mathbf{Claim~2}$: $U=\bigcup_{n\in\mathbb{N}}U_n$ and $V=\bigcup_{n\in\mathbb{N}}V_n$.
		
		It is easy to see that $\bigcup_{n\in\mathbb{N}}U_n\subseteq U$. Conversely, for any given $n\in \mathbb{N}$, $x_n\in W_k\cap U$ owing to $W_k\in \mathcal{U}_k$ for any $k\geq n$, which yields that $x_n\in U_n$. Therefore, $U=\bigcup_{n\in\mathbb{N}}U_n$. The result that $V=\bigcup_{n\in\mathbb{N}}V_n$ follows directly by applying a similar analysis.
		
		\textcolor{black}{We observe that} $(U_n,V_n)\leq (W_n,W_n')\in (L^*\times P^*)\backslash \mathcal{U}$. Then the fact that $(L^*\times P^*)\backslash \mathcal{U}$ is a lower set \textcolor{black}{implies} $(U_n,V_n)\in (L^*\times P^*)\backslash \mathcal{U}$. By the construction of $U_n$, $V_n$, we know that $(U_n,V_n)_{n\in \mathbb{N}}$ is a directed subset of $(L^*\times P^*)\backslash \mathcal{U}$. It follows that $\sup_{n\in \mathbb{N}}(U_n,V_n)=(U,V)\in (L^*\times P^*)\backslash \mathcal{U}$ from the Scott closedness of $(L^*\times P^*)\backslash \mathcal{U}$. This is a contradiction to the assumption that $(U,V)\in \mathcal{U}$.
		
		\textcolor{black}{In conclusion}, the product topology of $L^*\times P^*$ coincides with the Scott topology of the product poset $L^*\times P^*$.
	\end{proof}
\end{lemma}
\begin{theorem}\label{sober}
	For any countable poset $L$,  $\Sigma(\sigma(L))$ is a sober space.
	\begin{proof}
		Let $L^*=\sigma(L)$. From Lemma \ref{product}, we can obtain that $\Sigma(L^*\times L^*)=\Sigma L^*\times \Sigma L^*$. This means that $L^*$ is a sup semilattice such that the sup operation is jointly Scott-continuous by \cite[Corollary II-1.12]{gierz03}, as a result, $\Sigma L^*$ is sober.
	\end{proof}
\end{theorem}
\begin{theorem}
	For any two countable posets $L, P$, write $L^*=\sigma (L)$ and $P^*=\sigma(P)$, then the Scott space of  $L^*\times P^*$ is sober.
	\begin{proof}
		In the light of Theorem \ref{sober}, we know that $\Sigma L^*$ and $\Sigma P^*$ are both sober.  Due to \cite[Exercise 0-5.16.]{gierz03}, we get that products of sober spaces are sober. This reveals that the product of $\Sigma L^*\times \Sigma P^*$ is sober. It turns out that $\Sigma(L^*\times P^*)$ is sober via Lemma \ref{product}.
	\end{proof}
\end{theorem}
Kou's directed determined spaces \textcolor{black}{\cite{yukou15}} generalize dcpo's endowed with the Scott topology, where a $T_0$ space is said to be a \emph{directed determined space}, if an upper set $U$ in the specialization order is open iff \textcolor{black}{$\overline{D}\cap U\neq \emptyset$} implies that $D\cap U\neq \emptyset$ for any directed subset $D$ of $X$. Directed determined spaces have been found to be quite useful in constructing free objects in the category of dcpo's. Applying a similar discussion to directed determined spaces, the results in the remark below \textcolor{black}{follow} immediately.

\begin{remark}
	Let $X,Y$ be two countable directed determined spaces, $X^*=\mathcal{O}(X)$ and $Y^*=\mathcal{O}(Y)$. Then the following statements hold:
	\begin{enumerate}
		\item $\Sigma X^*\times \Sigma Y^*=\Sigma(X^*\times Y^*)$;
		
		\item $\Sigma X^*$ is sober;
		
		\item $\Sigma(X^*\times Y^*)$ is sober.
	\end{enumerate}
\end{remark}

\section{Two sober dcpo's with a non-sober dcpo product}

In this section, we give two dcpo's $P_1$ (Lemma~\ref{P1}) and $P_2$ (Lemma~\ref{P2}) that are sober in the Scott topology, and we will see that the Scott topology on their product poset $P_1\times P_2$ is not sober (Theorem~\ref{productisnotsober}). As mentioned in the Introduction, this solves a long-lasting question. Our construction of $P_1$ and $P_2$ go through a list of posets $\mathbb{N}^{<\mathbb{N}}$, $M$, $L$ (Example \ref{nandnn}) and poset $B$ (Example \ref{asdfdf}). 
In what follows, we use $\mathbb{N}$ to denote the set of all non-negative integers.

\begin{example}\label{nandnn}~
\begin{enumerate}

\item Let $\mathbb{N}^{<\mathbb{N}}$ be the poset  of all nonempty finite words (or, finite strings) over $\mathbb{N}$, with the {\sl prefix order} $\le$ defined as:
for $x=a_1a_2\cdots a_n$ and $ y= b_1b_2\cdots b_m$ in $\mathbb{N}^{<\mathbb{N}}$,
$x\le y \mbox{ if and only if } n\le m \mbox{ and } a_i =b_i \mbox{ for all } 1\le i\le n$. 
\item Let $M=\mathbb{N} \cup \mathbb{N}^{<\mathbb{N}}$. The order $\leq_M$ on $M$ is defined as $x\leq_M y$ if either $x\leq y$ in $\mathbb{N}^{<\mathbb{N}}$, or $x\leq y$ in $\mathbb{N}$ with respect to the usual order on $\mathbb N$. \textcolor{black}{See Figure~\ref{posetmp}}. 
\item Let $L=(\{(a, b)\in \mathbb{N}\times \mathbb{N}:a<b\}\times M)\cup \{\top\}$. \textcolor{black}{We define an order $\leq$ on $L$ as follows: $((a, b), x)\leq ((a', b'), y)$ if $a=a',b=b'$, $x\leq_M y$ in $M$, and $u\leq \top$ for all $u\in L$.} See Figure~\ref{posetlpic}.\\
We may write $((a, b), x)\in L$ simply as $x_{a,b}$ for $x\in M$  and $\{(a,b)\}\times M$ simply as $M_{a,b}$. 
\end{enumerate}
\end{example}

In order to define our poset $B$, we need to first fix an injective function~$i$ and a special function~$f$ described in the following remark.

\begin{remark}\label{i}

By Remark 2.6 of \cite{Miaoxl2023}, there is an injection $i:\{(a,b)\in \mathbb{N}\times \mathbb{N}:a<b\}\rightarrow \mathcal{P}(\mathbb{N})$ satisfying

\begin{enumerate}
\item $i(m_1, n_1)\cap i(m_2, n_2)=\emptyset$ if  $(m_1, n_1)\neq (m_2, n_2)$;
\item $n < k$ for every  $k\in i(m, n)$;
\item for each $(m,n)\in \mathbb{N}\times \mathbb{N}$ \textcolor{black}{and $m<n$}, there exists a monotone injective function $f_{m, n}: \mathbb{N}^{<\mathbb{N}}\to i(m, n)$. 
\end{enumerate}
\end{remark}

\begin{example}\label{asdfdf}~
Let $B=\mathbb{N}\times \mathbb{N}\times L$.  For each $s\in \mathbb{N}^{<\mathbb{N}}$ with length equaling 1, we regard it  as a natural number, also denoted by $s$. We now define four binary relations $\sqsubset_{1},  \sqsubset_{2}, \sqsubset_{3}$ and $\sqsubset_{4}$ on the set $B$ as follows:
\begin{itemize}
\item \textcolor{black}{$(m_1, n_1, x_{a_1,b_1})\sqsubset_{1}(m_2, n_2,y_{a_2,b_2})$} if \textcolor{black}{$m_1=m_2, n_1=n_2$ and $x_{a_1,b_1}< y_{a_2,b_2}$ holds in $L$. In other words, $(m, n, x)\sqsubset_{1} (m, n, y)$ if $x< y$ holds in $L$.}
\item  $(m_1, n_1,x_{a_1,b_1})\sqsubset_{2}(m_2,n_2,y_{a_2,b_2})$ if $y_{a_2,b_2}=\top$, $x\in \mathbb{N}^{<\mathbb{N}}$ and $m_1=a_1,n_2=n_1+1, m_2=f_{a_1,b_1}(x)$.
In other words, $(a, n, x_{a, b})\sqsubset_{2} (f_{a, b}(x), n+1,\top)$ for any \textcolor{black}{$x\in \mathbb{N}^{<\mathbb{N}}$}.
\item $(m_1,n_1,x_{a_1,b_1})\sqsubset_{3}(m_2,n_2,y_{a_2,b_2})$ if $y_{a_2, b_2}=\top$, $x\in \mathbb{N}$ and $m_1=b_1, n_2=n_1+1, m_2=f_{a_1,b_1}(x)$. In other words,
$(b, n, x_{a, b})\sqsubset_{3}(f_{a, b}(x),n+1,\top)$ for any $x\in \mathbb{N}$. {Here $x\in \mathbb{N}$, and in the definition of $f_{a_1, b_1}(x)$, $x\in\mathbb{N}$ is taken as an element of $\mathbb{N}^{<\mathbb{N}}$ with length equaling 1.}
\item $(m_1,n_1, x_{a_1,b_1})\sqsubset_{4}(m_2,n_2,y_{a_2,b_2})$ if $y_{a_2,b_2}=\top$, $x\in \mathbb{N}$ and there exists $s\in  \mathbb{N}^{<\mathbb{N}}$ such that $m_1=f_{a_1,b_1}(s),m_2=f_{a_1,b_1}(s.x), n_1=n_2$. In other words, $(f_{a, b}(s), n, x_{a, b})\sqsubset_{4}(f_{a, b}(s.x), n, \top)$ for any $s\in  \mathbb{N}^{<\mathbb{N}}$ and $x\in \mathbb{N}$.
\end{itemize}
Now let
\begin{itemize}
\item $\sqsubset=\sqsubset_1\cup \sqsubset_2\cup\sqsubset_3\cup\sqsubset_4\cup\sqsubset_1;\sqsubset_2\cup\sqsubset_1;\sqsubset_3\cup \sqsubset_1;\sqsubset_4$,
and then let $\sqsubseteq$ be  $\sqsubset\cup =$. 
\end{itemize}
Here $=$ denotes the identity relation on $B$, and the relation  $\sqsubset_1;\sqsubset_2$ is the composition of  $\sqsubset_1$ and $\sqsubset_2$:
that is,  $(m_1, n_1, x_{a_1, b_1})\sqsubset_{1} (m_2, n_2, x_{a_2, b_2}) \sqsubset_{2} (m_3, n_3, x_{a_3, b_3})$ would imply that $(m_1, n_1, x_{a_1, b_1})\sqsubset_1;\sqsubset_2  (m_3, n_3, x_{a_3, b_3})$.  It is clear that $\sqsubseteq$ is a partial order on $B$.
\end{example}

\begin{figure}[t]
	\centering
	\begin{tikzpicture}[scale=1]
		\foreach \x in {0,1,2} \fill (0,\x) circle (2pt);
		\fill (0,4) circle (0pt);
		\draw (0,0)--(0,2);
		\draw[densely dashed] (0,2)--(0,4);
		\fill (2.5,0) circle (2pt);
		\fill (1.5,1) circle (2pt);
		\fill (3.5,1) circle (2pt);
		\foreach \x in {1,2,3,4} \fill (\x,2) circle (2pt);
		\draw (2.5,0)--(1.5,1);
		\draw (2.5,0)--(3.5,1);
		\draw (1.5,1)--(1,2);
		\draw (1.5,1)--(2,2);
		\draw (3.5,1)--(3,2);
		\draw (3.5,1)--(4,2);
		\fill (3,0) circle (0pt);
		\fill (4.8,0) circle (0pt);
		\draw[densely dashed] (3,0)--(4.8,0);
		\fill (3.8,1) circle (0pt);
		\fill (4.8,1) circle (0pt);
		\draw[densely dashed] (3.8,1)--(4.8,1);
		\fill (2.2,2) circle (0pt);
		\fill (2.8,2) circle (0pt);
		\draw[densely dashed] (2.2,2)--(2.8,2);
		\fill (4.2,2) circle (0pt);
		\fill (4.8,2) circle (0pt);
		\draw[densely dashed] (4.2,2)--(4.8,2);
		\fill (1.5,2) circle (0pt);
		\fill (1.5,4) circle (0pt);
		\draw[densely dashed] (1.5,2)--(1.5,4);
		\fill (3.5,2) circle (0pt);
		\fill (3.5,4) circle (0pt);
		\draw[densely dashed] (3.5,2)--(3.5,4);
		\node (1) at (0,0) [left] {\small  $0$};
		\node (2) at (0,1) [left] {\small  $1$};
		\node (3) at (0,2) [left] {\small  $2$};
		\node (4) at (2.5,0) [left] {\small  $0$};
		\node (5) at (1.5,1) [left] {\small  $00$};
		\node (6) at (3.5,1) [left] {\small  $01$};
		\node (7) at (1,2) [above] {\small  $000$};
		\node (8) at (2,2) [above] {\small  $001$};
		\node (9) at (3,2) [above] {\small  $010$};
		\node (10) at (4,2) [above] {\small  $011$};
	\end{tikzpicture}
	\caption{The poset $M$ \label{posetmp}}
\end{figure}

\begin{figure}[t]
	\centering
	\begin{tikzpicture}[scale=0.8]
		\foreach \x in {0,1,2} \fill (0,\x) circle (2pt);
		\fill (0,4) circle (0pt);
		\draw (0,0)--(0,2);
		\draw[densely dashed] (0,2)--(9.4,5);
		\fill (2.5,0) circle (2pt);
		\fill (1.5,1) circle (2pt);
		\fill (3.5,1) circle (2pt);
		\foreach \x in {1,2,3,4} \fill (\x,2) circle (2pt);
		\draw (2.5,0)--(1.5,1);
		\draw (2.5,0)--(3.5,1);
		\draw (1.5,1)--(1,2);
		\draw (1.5,1)--(2,2);
		\draw (3.5,1)--(3,2);
		\draw (3.5,1)--(4,2);
		\fill (3,0) circle (0pt);
		\fill (4.8,0) circle (0pt);
		\draw[densely dashed] (3,0)--(4.8,0);
		\fill (3.8,1) circle (0pt);
		\fill (4.8,1) circle (0pt);
		\draw[densely dashed] (3.8,1)--(4.8,1);
		\fill (2.2,2) circle (0pt);
		\fill (2.8,2) circle (0pt);
		\draw[densely dashed] (2.2,2)--(2.8,2);
		\fill (4.2,2) circle (0pt);
		\fill (4.8,2) circle (0pt);
		\draw[densely dashed] (4.2,2)--(4.8,2);
		\fill (1.5,2) circle (0pt);
		\fill (1.5,4) circle (0pt);
		\draw[densely dashed] (1.5,2)--(9.4,5);
		\fill (3.5,2) circle (0pt);
		\fill (3.5,4) circle (0pt);
		\draw[densely dashed] (3.5,2)--(9.4,5);
		\node (1) at (0,0) [left] {\small \tiny$0_{0,1}$};
		\node (2) at (0,1) [left] {\small \tiny$1_{0,1}$};
		\node (3) at (0,2) [left] {\small  \tiny$2_{0,1}$};
		\node (4) at (2.5,0) [left] {\small   \tiny$0_{0,1}$};
		\node (5) at (1.5,1) [left] {\small   \tiny$00_{0,1}$};
		\node (6) at (3.5,1) [left] {\small   \tiny$01_{0,1}$};
		\node (7) at (1,2) [above] {\small   \tiny$000_{0,1}$};
		\node (8) at (2,2) [above] {\small   \tiny$001_{0,1}$};
		\node (9) at (3,2) [above] {\small   \tiny$010_{0,1}$};
		\node (10) at (4,2) [above] {\small   \tiny$011_{0,1}$};
		\node (11) at (2.5,0) [below] {$\{(0,1)\}\times M$};
		\foreach \x in {0,1,2} \fill (6,\x) circle (2pt);
		\fill (6,4) circle (0pt);
		\draw (6,0)--(6,2);
		\draw[densely dashed] (6,2)--(9.4,5);
		\fill (8.5,0) circle (2pt);
		\fill (7.5,1) circle (2pt);
		\fill (9.5,1) circle (2pt);
		\foreach \x in {7,8,9,10} \fill (\x,2) circle (2pt);
		\draw (8.5,0)--(7.5,1);
		\draw (8.5,0)--(9.5,1);
		\draw (7.5,1)--(7,2);
		\draw (7.5,1)--(8,2);
		\draw (9.5,1)--(9,2);
		\draw (9.5,1)--(10,2);
		\fill (9,0) circle (0pt);
		\fill (10.8,0) circle (0pt);
		\draw[densely dashed] (9,0)--(10.8,0);
		\fill (9.8,1) circle (0pt);
		\fill (10.8,1) circle (0pt);
		\draw[densely dashed] (9.8,1)--(10.8,1);
		\fill (8.2,2) circle (0pt);
		\fill (8.8,2) circle (0pt);
		\draw[densely dashed] (8.2,2)--(8.8,2);
		\fill (10.2,2) circle (0pt);
		\fill (10.8,2) circle (0pt);
		\draw[densely dashed] (10.2,2)--(10.8,2);
		\fill (7.5,2) circle (0pt);
		\fill (7.5,4) circle (0pt);
		\draw[densely dashed] (7.5,2)--(9.4,5);
		\fill (9.5,2) circle (0pt);
		\fill (9.5,4) circle (0pt);
		\draw[densely dashed] (9.5,2)--(9.4,5);
		\node (1) at (6,0) [left] {\small \tiny$0_{0,2}$};
		\node (2) at (6,1) [left] {\small \tiny$1_{0,2}$};
		\node (3) at (6,2) [left] {\small  \tiny$2_{0,2}$};
		\node (4) at (8.5,0) [left] {\small   \tiny$0_{0,2}$};
		\node (5) at (7.5,1) [left] {\small   \tiny$00_{0,2}$};
		\node (6) at (9.5,1) [left] {\small   \tiny$01_{0,2}$};
		\node (7) at (7,2) [above] {\small   \tiny$000_{0,2}$};
		\node (8) at (8,2) [above] {\small   \tiny$001_{0,2}$};
		\node (9) at (9,2) [above] {\small   \tiny$010_{0,2}$};
		\node (10) at (10,2) [above] {\small   \tiny$011_{0,2}$};
		\node (11) at (8.5,0) [below] {$\{(0,2)\}\times M$};
		\fill (11.5,0) circle (0pt);
		\fill (13,0) circle (0pt);
		\draw[densely dashed] (11.5,0)--(13,0);
		\fill (11.5,1) circle (0pt);
		\fill (13,1) circle (0pt);
		\draw[densely dashed] (11.5,1)--(13,1);
		\fill (11.5,2) circle (0pt);
		\fill (13,2) circle (0pt);
		\draw[densely dashed] (11.5,2)--(13,2);
		\fill (11.5,3) circle (0pt);
		\fill (13,3) circle (0pt);
		\foreach \x in {0,1,2} \fill (14,\x) circle (2pt);
		\fill (14,4) circle (0pt);
		\draw (14,0)--(14,2);
		\draw[densely dashed] (14,2)--(9.4,5);
		\fill (16.5,0) circle (2pt);
		\fill (15.5,1) circle (2pt);
		\fill (17.5,1) circle (2pt);
		\foreach \x in {15,16,17,18} \fill (\x,2) circle (2pt);
		\draw (16.5,0)--(15.5,1);
		\draw (16.5,0)--(17.5,1);
		\draw (15.5,1)--(15,2);
		\draw (15.5,1)--(16,2);
		\draw (17.5,1)--(17,2);
		\draw (17.5,1)--(18,2);
		\fill (17,0) circle (0pt);
		\fill (18.8,0) circle (0pt);
		\draw[densely dashed] (17,0)--(18.8,0);
		\fill (17.8,1) circle (0pt);
		\fill (17.8,1) circle (0pt);
		\draw[densely dashed] (17.8,1)--(18.8,1);
		\fill (16.2,2) circle (0pt);
		\fill (16.8,2) circle (0pt);
		\draw[densely dashed] (16.2,2)--(16.8,2);
		\fill (18.2,2) circle (0pt);
		\fill (18.8,2) circle (0pt);
		\draw[densely dashed] (18.2,2)--(18.8,2);
		\fill (15.5,2) circle (0pt);
		\fill (15.5,4) circle (0pt);
		\draw[densely dashed] (15.5,2)--(9.4,5);
		\fill (17.5,2) circle (0pt);
		\fill (17.5,4) circle (0pt);
		\draw[densely dashed] (17.5,2)--(9.4,5);
		\node (1) at (14,0) [left] {\small \tiny$0_{a,b}$};
		\node (2) at (14,1) [left] {\small \tiny$1_{a,b}$};
		\node (3) at (14,2) [left] {\small  \tiny$2_{a,b}$};
		\node (4) at (16.5,0) [left] {\small   \tiny$0_{a,b}$};
		\node (5) at (15.5,1) [left] {\small   \tiny$00_{a,b}$};
		\node (6) at (17.5,1) [left] {\small   \tiny$01_{a,b}$};
		\node (7) at (15,2) [above] {\small   \tiny$000_{a,b}$};
		\node (8) at (16,2) [above] {\small   \tiny$001_{a,b}$};
		\node (9) at (17,2) [above] {\small   \tiny$010_{a,b}$};
		\node (10) at (18,2) [above] {\small   \tiny$011_{a,b}$};
		\node (11) at (16.5,0) [below] {$\{(a,b)\}\times M$};
		\fill (9.4,5) circle (2pt);
		\node (1) at (9.4,5) [above] {\small $\top$};
		
	\end{tikzpicture}
	\caption{The poset $L$ \label{posetlpic}}
	\end{figure}

\input{figure3}

\input{figure4}

\input{figure5}

 \textcolor{black}{For a clearer understanding of the sobriety of $P_1$ and $P_2$ that we will define below, we present Figures~\ref{b-leq-2}, \ref{b-leq-3} and \ref{b-leq-4}, depicting the relations $\sqsubset_2, \sqsubset_3$ and $\sqsubset_4$ on $B$, respectively.}

\input{figure6}

\input{figure7}

\begin{remark}\label{directed}
By the definition of  the order $\sqsubseteq$ on $B$, it follows easily that  a non-trivial directed subset $D$ of $B$  is of the form
$$D=\{(m_0, n_0, (x_k)_{a_0, b_0}): k\in N\}$$
for some fixed $m_0, n_0, a_0, b_0\in \mathbb{N}$, where either $(x_k)_{k\in N}$ is a cofinal subset of $\mathbb{N}$ or a non-trivial directed subset in $\mathbb{N}^{<\mathbb{N}}$. In either case, the indexed set $N$ is a cofinal subset of $\mathbb{N}$.
Hence $\bigvee D=(m_0, n_0, \top)$ always exists in $B$.
\end{remark}

In the next example, we will define our dcpo $P_1$. 
\begin{example}
Let first $\mathbb{N}^{\mathbb{N}}$ denote the set of all mappings from $\mathbb{N}$ to $\mathbb{N}$. Define the  order  $\le$ on $\mathbb{N}^{\mathbb{N}}\times \mathbb{N} $  as: for any $(f,n),(g,m)$ in $\mathbb{N}^{\mathbb{N}}\times \mathbb{N}$,
$$(f,n)\le (g,m) \mbox{ if and only if } f = g \mbox{ and } n\le m.$$
Now let $P_1=(\mathbb{N}^{\mathbb{N}}\times \mathbb{N})\cup B\cup \{\top_1\}$. 

We define the partial order $\leq$ on $P_1$ as the one generated by relations $<_{1},  <_{2}, <_{3}$ and $<_{4}$  as follows:
\begin{itemize}
\item $<_{1}=\le \setminus =$, where $\le$ is the order relation  on $\mathbb{N}^{\mathbb{N}}\times \mathbb{N}$.
\item $<_{2}=\sqsubseteq \setminus =$ , where  $\sqsubseteq$ is   the order relation on $B$.
\textcolor{black}{\item $x<_{3}y$ if $x=(f, n)\in \mathbb{N}^{\mathbb{N}}\times \mathbb{N}$, $y=(f(n),n,\top)\in B$. We could see that $(f, n)<_{3}(f(n), n,\top)$.}
\textcolor{black}{\item $x<_{4}y$ if $x\in P_1\backslash\{\top_1\}$, $y=\top_1$.}
\end{itemize}
Let
$$< = <_{1}\cup <_{2}\cup <_{3}\cup <_{4}\cup <_{1};<_{3}, $$
and then the order $\leq$ on $P_1$ is defined as $ < \cup =$,  the union of $ <$ and the identity relation. \textcolor{black}{(See Figure \ref{FigureP1})}
\end{example}

We continue to define our dcpo~$P_2$. 

\begin{example}

Now take one bijection $\phi: \mathbb{N}\rightarrow \{(a,b)\in \mathbb{N}\times \mathbb{N}:a<b\}$. This is possible because $\{(a, b)\in \mathbb{N}\times\mathbb{N}: a < b\}$ is countable. For each $n\in \mathbb{N}$, let $E_n=i(\phi(n))$, that is if $\phi(n)=(a, b)$, then $E_n=i(a, b)$. For clearness, we may also use $E_{n}^{a, b}$ for $E_n$. Let $(\bigcup_{n\in \mathbb{N}}E_{n})^{\mathbb{N}}$ denote all mappings from $\mathbb{N}$ to $\bigcup_{n\in \mathbb{N}}E_{n}$. Set $X=\{f\in(\bigcup_{n\in \mathbb{N}}E_{n})^{\mathbb{N}}:\forall n\in \mathbb{N}, f(n)\in E_{n}\}$. Define the order $\le$ on $(X\times \mathbb{N})\times \mathbb{N}$  as: for any $((f,n),k),((g,m),p)$ in $(X\times \mathbb{N})\times \mathbb{N}$, 
\begin{itemize}
\item $((f,n),k)\le ((g,m),p) \mbox{ if and only if } n\le m \mbox{ and } f=g,k=p.$
\end{itemize}
We will write $((f,n),k)$ as $(f,n,k)$ for any $((f,n),k)\in (X\times \mathbb{N})\times \mathbb{N}$. 

Finally, let $P_2=((X\times \mathbb{N})\times \mathbb{N})\cup B\cup \{\top_2\}$. We define the relations $<_{1},  <_{2}, <_{3}$ and $<_{4}$ on $P_2$ as follows:
\begin{itemize}
\item $<_{1} $ is the  order relation, excluding the equality, on $(X\times \mathbb{N})\times \mathbb{N}$.
\item $<_{2} = \sqsubset$ on $B$.
\item $ x<_{3}y $ if $x=(f,n,k)\in(X\times \mathbb{N})\times \mathbb{N}$, $y=(f(n),k,\top)\in B$. In other words, $(f,n,k)<_{3}(f(n),k,\top)$.
\item $x<_{4}y$ if $x\in P_2\backslash\{\top_2\}$, $y=\top_2$.
\end{itemize}
Let $<=<_{1}\cup <_{2}\cup <_{3}\cup <_{4}\cup <_{1};<_{3}$ and $\leq = <\cup =$, the union of $ < $ and the \textcolor{black}{identity} relation $=$. This is the order that we impose on $P_2$. \textcolor{black}{(See Figure \ref{FigureP2})}
\end{example}

\input{figure8}

\input{figure9}

From the definitions of $P_1$ and $P_2$ and Remark~\ref{directed}, \textcolor{black}{we can
easily prove} the following proposition.
\textcolor{black}{
\begin{proposition}\label{dcpo}~
\begin{enumerate}
    \item  Every directed subset $D$ of $P_1$ without greatest elements, is either of the form $\{(f_0,n): n\in N\}$ (where $N$ is an infinite subset of $\mathbb{N}$), or  entirely contained in $B\subseteq P_1$.
    \item Every directed subsets $D$ of $P_2$ without greatest elements, is either of the form $\{(f_0,n,k_0): {n\in N}\}$ (where $N$ is an infinite subset of $\mathbb{N}$), or  entirely contained in $B\subseteq P_2$.
\end{enumerate}
Hence, the posets $P_1$ and $P_2$ are dcpo's.
\end{proposition}}

Note that for any non-empty lower subset $A$ of a dcpo $P$, $A$ is Scott closed iff for any directed subset $D\subseteq A$ without maximal elements, \textcolor{black}{$\sup D$ belongs to $A$}. Furthermore, \textcolor{black}{for any Scott closed} set $A$ of a dcpo $P$, $A=\da \max A$, where $\max A$ is the set of all maximal points of $A$. Also, if $A$ is irreducible and  not the closure of a singleton, $\max A$ must be an infinite set.

\begin{lemma}\label{P1}
$\Sigma P_1$ is sober.
\begin{proof}
Suppose that  $\Sigma P_1$ is not sober. Then there exists an irreducible closed subset $A$ of $P_1$ such that $A\notin \{\da x: x\in P_1\}$. By the previous remark,  $A=\da \max A$ and $\max A$ is infinite.

\textbf{Claim $1$}: For any $x\in \max A$, there exists a directed subset $D\subseteq B$ without maximal elements, such that $x=\sup D$ and $x\in \max B$. Thus $\max A\subseteq \max B$.

To see that, first note that $A=\da \max A=\da x\cup cl(\max A\backslash\{x\})$. Then, as $A$ is irreducible,  $A=\da x$ or $A=cl(\max A\backslash\{x\})$. By the assumption, $A\not=\da x$, we have
 $A=cl(\max A\backslash\{x\})=cl(\da(\max A\backslash\{x\}))$.
 For any non-trivial directed subset $D$ of $\da (\max A\backslash\{x\})$, we know that $C=D\cap (P_1\backslash B)$ is not a cofinal subset of $D$ (otherwise,  $\top_1=\bigvee C =\bigvee D\in A$). Then $D\cap B$ is a cofinal subset of $D$. By the definition of the order on $B$, we have that $\sup D\in \max B$. \textcolor{black}{It follows}  that $\sup D\in \max A$. Now, because  $\da (\max A\backslash\{x\})$ is not Scott closed, we can find a non-trivial directed subset $D$ of $\da (\max A\backslash\{x\})$ such that $\sup D\notin \da(\max A\backslash \{x\})$. By the above argument, $\sup D\in \max A$, hence  $\sup D=x\in \max B$.

Define $F=\{n\in \mathbb{N}:\max A\cap L_n\neq \emptyset\}$, where $L_n=\mathbb{N}\times \{n\}\times L$.

\textbf{Claim $2$}: $F$ is an infinite subset of $\mathbb{N}$.

Assume $F$ is a finite set. Let  $\max F=n^*$. It follows that

$$A= A\cap (\bigcup_{k\leq n^*-1}\da L_k)\cup cl(\max A\cap L_{n^*}).$$

\textbf{Claim $2.1$}: $A\cap (\bigcup_{k\leq n^*-1}\da L_k)$ is a Scott closed subset of $P_1$.

Clearly,  $A\cap (\bigcup_{k\leq n^*-1}\da L_k)$ is a lower set. Let $D$ be a directed subset of $A\cap (\bigcup_{k\leq n^{*}-1}\da L_k)$ without maximal elements. As $\top_1$ does not belong to $A$, $\sup D\neq \top_1$. We know that $D\cap(P_1\backslash B)$ is not a cofinal subset of $D$.
This infers that $D\cap B$ is a cofinal subset of $D$. By Remark \ref{directed},  $D\cap B =\{ (m_0,n_0,(x_k)_{a_0,b_0}):{k\in N}\}$ for some fixed $m_0, n_0,a_0,b_0\in \mathbb{N}$, where $(x_k)_{k\in N}$ is either a cofinal subset of $\mathbb{N}$ or a non-trivial directed subset in $\mathbb{N}^{<\mathbb{N}}$.

For each $k\in\mathbb{N}$,  $(m_0,n_0, (x_k)_{a_0,b_0})\in \da L_{n_k}=\da \max L_{n_k}$ for some $n_k\leq n^*-1$, and it means that $(m_0,n_0,(x_k)_{a_0,b_0})<_2 (m_k,n_k,\top)$ in $P_1$ for some $m_k\in \mathbb{N}$. Now it remains to show that $\sup D=\sup_{k\in N}(m_0,n_0,(x_k)_{a_0,b_0})=(m_0,n_0,\top)\in A\cap (\bigcup_{k\leq n^*-1}\da L_k)$.
We will only consider the non-trivial case when the set $\{(m_k,n_k,\top): k\in \mathbb{N}\}$ is  infinite and $(m_0,n_0,\top)$ does not belong to $\{(m_k,n_k,\top): k\in \mathbb{N}\}$. Without loss of generality, we assume that $(m_k, n_k, \top)\not=(m_h, n_h, \top)$ when $k\not=h$.
By the order of $P_1$ and $(m_0,n_0)\neq (m_k,n_k)$, $<_2\subseteq \sqsubset_2\cup\sqsubset_3\cup\sqsubset_4\cup\sqsubset_1;\sqsubset_2\cup\sqsubset_1;\sqsubset_3\cup \sqsubset_1;\sqsubset_4$ in this case.
We let

\begin{center}
	$A_{\sqsubset_{r}} = \{k\in N:(m_0, n_0,(x_k)_{a_0,b_0})\sqsubset_r (m_k, n_k,\top)\}$,
\end{center}
where $\sqsubset_{r}\in \{\ \sqsubset_{2},\ \sqsubset_{3},\ \sqsubset_{4},\ \sqsubset_{1};\sqsubset_{2},\ \sqsubset_{1};\sqsubset_{3},\ \sqsubset_{1};\sqsubset_{4}\}$. It follows that at least one of the six relations is infinite because that $N$ is infinite. 

{Case $1$: $A_{\sqsubset_2}$ or $A_{\sqsubset_{1};\sqsubset_{2}}$ is infinite. Then $n_k=n_0+1\leq n^*-1$ for each $k$ in $A_{\sqsubset_{2}}$ or $A_{\sqsubset_{1};\sqsubset_{2}}$. It follows that $n_0\leq n^*-1$. Hence $\sup D=\sup_{k\in A_{\sqsubset_2}}(m_0,n_0,(x_k)_{a_0,b_0})=(m_0,n_0,\top)$, or $\sup D=\sup_{k\in A_{\sqsubset_1;\sqsubset_2}}(m_0,n_0,(x_k)_{a_0,b_0})=(m_0,n_0,\top)$. Since $D$ is a directed subset contained in the Scott closed set~$A$, $\sup D=(m_0,n_0,\top)\in A$.
	Note that $(m_0,n_0,\top)\in A\cap L_{n_0}\subseteq A\cap (\bigcup_{k\leq n^*-1}\da L_k)$. Hence, $\sup D\in A\cap(\bigcup_{k\leq n^*-1}\da L_k)$.}

{Case $2$: $A_{\sqsubset_3}$ or $A_{\sqsubset_{1};\sqsubset_{3}}$ is infinite. Then $n_k=n_0+1\leq n^*-1$ for each $k\in A_{\sqsubset_{3}}$ or $k\in A_{\sqsubset_{1};\sqsubset_{3}}$. The remaining part of the  proof is similar to that in Case $1$.}

{Case $3$: $A_{\sqsubset_4}$ or $A_{\sqsubset_{1};\sqsubset_{4}}$ is infinite. Then $n_k=n_0\leq n^*-1$ for each $k\in A_{\sqsubset_{4}}$ or $k\in A_{\sqsubset_{1};\sqsubset_{4}}$. The remaining part of the  proof is similar to that in Case $1$.}

So we have  $A\cap (\bigcup_{k\leq n^*-1}\da L_k)$ is a Scott closed subset of $P_1$. Remember that $A$ is irreducible, so we can get that $A =  A\cap (\bigcup_{k\leq n^*-1}\da L_k)$ or $A = cl(\max A\cap L_{n^*})$. Since $\max A\cap L_{n^*}\neq \emptyset$,  we must have that  $A=cl(\max A\cap L_{n^*})$.

Choose an element $(m_1,n^*,\top)\in \max A\cap L_{n^*}$. By  Claim $1$,  there is a non-trivial directed subset $D$ of $B\cap \da(\max A\backslash\{(m_1,n^*,\top)\})$ such that $\sup D=(m_1,n^*,\top)$.
By Remark \ref{directed}, we have $D=\{(m_1,n^{*},(y_k)_{a_1,b_1}):{k\in N}\}$, where $N$ is an infinite subset of $\mathbb{N}$. For each $k\in N$, there exists $(a_k,n_k,\top)\in \max A\backslash\{(m_1,n^*,\top)\}$ such that  $(m_1,n^{*},(y_k)_{a_1,b_1})\leq (a_k,n_k,\top)$. This leads to that $n^*\leq n_k$ with a similar analysis as in Claim 2.1. Note that $n_k\in F$. Then $n_k\leq \max F=n^*$. Thus $n_k=n^*$ for all $k\in N$.
It follows that the set $M=\{(a_k,n^*,\top):k\in N\}\subseteq \max A\backslash\{(m_1,n^*,\top)\}$ is infinite. To see this, assume on the contrary that $M$ is finite. Then $D\subseteq \da M$, and $\da M$ is Scott closed, which yields that $\sup D=(m_1,n^*,\top)\in \da M$. This means that there is $(a_k,n^*,\top)\in M$ such that $(m_1,n^*,\top)\leq(a_k,n^*,\top)$. Hence we know that $(m_1,n^*,\top)=(a_k,n^*,\top)$ by the order on $P_1$, and this contradicts that $(a_k,n_k,\top)\in \max A\backslash\{(m_1,n^*,\top)\}$. Now, 
from the definition of $<_2$ in $P_1$, we can deduce that $(m_1,n^{*},(y_k)_{a_1,b_1})<_2 (a_k,n^*,\top)$, and since the two elements share the same second component $n^*$, we know that $<_2~\in \{\sqsubset_4,\ \sqsubset_1;\sqsubset_4\}$. This means that there is $s_0\in \mathbb{N}^{<\mathbb{N}}$ satisfying that $m_1=f_{a_1,b_1}(s_0)$ and $a_k=f_{a_1,b_1}(s_0.z_k)$ for some $z_k\geq y_k$ in $\mathbb{N}$. Since $M$ is an infinite set, $\{a_k:k\in N\}$ is also infinite. It follows that there are two distinct  $k_1,k_2$ in $N$ such that  $z_{k_1}\neq z_{k_2}$ as $f_{a_1,b_1}$ is injective.

We set
\begin{center}
$A_{k_1}=\da \{(f_{a_1,b_1}(s),n^*,\top)\in \max A:s\geq s_0.z_{k_1}\}\cup (A\cap (\bigcup_{k\leq n^*-1}\da L_k)$);\\
$A_{k_2}=\da \{(f_{a_1,b_1}(s),n^*,\top)\in \max A:s\geq s_0.z_{k_2}\}\cup (A\cap (\bigcup_{k\leq n^*-1}\da L_k)$).
\end{center}

\textbf{Claim $2.2$}: $A_{k_1}$ is Scott closed in $P_1$.

Obviously, $A_{k_1}$ is a lower set of $P_1$. Let $D$ be a non-trivial directed subset of $A_{k_1}$. We consider  the following two cases for $D$.

Case $1$: $D\cap  A\cap (\bigcup_{k\leq n^*-1}\da L_k)$ is cofinal in $D$. By Claim $2.1$, $A\cap (\bigcup_{k\leq n^*-1}\da L_k)$ is Scott closed. Hence,
$$\sup D=\sup (D\cap  A\cap (\bigcup_{k\leq n^*-1}\da L_k))\in A\cap (\bigcup_{k\leq n^*-1}\da L_k) \subseteq A_{k_1}.$$

Case $2$: $D\cap (\da \{(f_{a_1,b_1}(s),n^*,\top)\in \max A:s\geq s_0.z_{k_1}\})$ is cofinal in $D$. Note that $\top_1\not\in A$, so the set $D\cap (P_1\backslash B)\cap (\da \{(f_{a_1,b_1}(s),n^*,\top)\in \max A:s\geq s_0.z_{k_1}\})$ is not cofinal in $D$. Then we know $D\cap B\cap (\da \{(f_{a_1,b_1}(s),n^*,\top)\in \max A:s\geq s_0.z_{k_1}\})$ is cofinal in $D$.  By  Remark \ref{directed}, we have $$ D\cap B\cap (\da \{(f_{a_1,b_1}(s),n^*,\top)\in \max A:s\geq s_0.z_{k_1}\})=\{(m_0, n_0,(x_k)_{a_0,b_0}):{k\in N}\}$$ for some fixed $m_0, n_0,a_0,b_0\in \mathbb{N}$, where $(x_k)_{k\in N}$ is a cofinal subset of $\mathbb{N}$ or a non-trivial directed subset in $\mathbb{N}^{<\mathbb{N}}$.
Then for each $k\in N$, there is $(m_k,n^{*},\top)\in \{(f_{a_1,b_1}(s), n^*,\top)\in \max A:s\geq s_0.z_{k_1}\}$ such that $(m_0,n_0, (x_k)_{a_0,b_0})\leq (m_k,n^*,\top)$, and with a similar analysis as in Claim~2.1, we know that $n_0\leq n^*$.

If $n_0\leq n^*-1$, then $\sup D=(m_0,n_0,\top)\in \max A\cap L_{n_0}\subseteq A\cap (\bigcup_{k\leq n^*-1}\da L_k)\subseteq A_{k_1}$.

Now let $n_0=n^*$. 
Then for each $k\in N$, there is $(m_k,n^{*},\top)\in \{(f_{a_1,b_1}(s), n^*,\top)\in \max A:s\geq s_0.z_{k_1}\}$ such that $(m_0,n_0, (x_k)_{a_0,b_0})\leq (m_k,n^*,\top)$. 
We claim $\sup D=\sup_{k\in N}(m_0,n_{0},(x_k)_{a_0,b_0})$ $=(m_0,n_{0},\top)\in A_{k_1}$. Again, we only consider the non-trivial case: $\{(m_k,n^*,\top):k\in N\}$ is infinite and $(m_0,n_{0},\top)\not\in \{(m_k,n^*,\top):k\in N\}$. By the definition of $<_2$ in $P_1$, we must have that $(m_0,n_{0},(x_k)_{a_0,b_0})<_2 (m_k,n^*,\top)$ holds. 
Again, the fact that $n_0=n^*$ implies that  $<_2\in \{\sqsubset_4,\ \sqsubset_1;\sqsubset_4\}$. It follows that we can find $s_1\in \mathbb{N}^{<\mathbb{N}}$ satisfying that $m_0=f_{a_0,b_0}(s_1)$ and $m_k=f_{a_0,b_0}(s_1.b_k)$ for some $b_k\geq x_k$ in $\mathbb{N}$, and $(m_k,n^{*},\top)\in \{(f_{a_1,b_1}(s), n^*,\top)\in \max A:s\geq s_0.z_{k_1}\}$ implies that $f_{a_0,b_0}(s_1.b_k)=m_k=f_{a_1,b_1}(s)$ with $s\geq s_0.z_{k_1}$ in $\mathbb{N}^{<\mathbb{N}}$.
This indicates that $(a_0,b_0)=(a_1,b_1)$ (note that if $(a_0, b_0)\not=(a_1, b_1)$, then $i(a_0, b_0)\cap i(a_1, b_1)=\emptyset$) and $s_1.b_k=s\geq s_0.z_{k_1}$. 
We know that $\{b_k:k\in N\}$ is infinite as $f_{a_0,b_0}$ is injective and we assumed that $\{(m_k,n^*,\top):k\in N\}$ is infinite.
Hence there are at least two different $b_k, b_{k'} $ in $\mathbb{N}$ such that $s_1.b_k\geq s_0.z_{k_1}$ and $s_1.b_{k'}\geq s_0.z_{k_1}$, so we can obtain $s_0.z_{k_1}\leq s_1$ by the order on $\mathbb{N}^{<\mathbb{N}}$; and therefore,
$\sup D=(m_0,n_0,\top)=(f_{a_1,b_1}(s_1), n^*,\top)\in \{(f_{a_1,b_1}(s), n^*,\top)\in \max A:s\geq s_0.z_{k_1}\}\subseteq A_{k_1}.$

Employing a similar argument, we can show that $A_{k_2}$ is also a Scott closed subset of $P_1$.  We also know that  $A =  A_{k_1}\cup A_{k_2}\cup \da(\max A\backslash (A_{k_1}\cup A_{k_2}))\cup (A\cap (\bigcup_{k\leq n^*-1}\da L_k)).$

\textbf{Claim $2.3$}: $\da(\max A\backslash (A_{k_1}\cup A_{k_2}))\cup (A\cap (\bigcup_{k\leq n^*-1}\da L_k))$ is Scott closed.

Let $D$ be a non-trivial directed subset of the lower set $\da(\max A\backslash (A_{k_1}\cup A_{k_2}))\cup (A\cap (\bigcup_{k\leq n^*-1}\da L_k))$. If $D\cap (A\cap (\bigcup_{k\leq n^*-1}\da L_k))$ is cofinal in $D$, then $\sup D=\sup (D\cap (A\cap (\bigcup_{k\leq n^*-1}\da L_k)))\in A\cap (\bigcup_{k\leq n^*-1}\da L_k))$ because $A\cap (\bigcup_{k\leq n^*-1}\da L_k)$ is Scott closed. Hence we will only need to discuss the case when $D\cap \da(\max A\backslash (A_{k_1}\cup A_{k_2}))$ is cofinal in $D$, or without loss of generality, when
$D\subseteq \da(\max A\backslash (A_{k_1}\cup A_{k_2}))$.

Note that $D\subseteq \da(\max A\backslash (A_{k_1}\cup A_{k_2}))\subseteq A$ and $\top_1\not\in A$. If $D\cap P_1\backslash B$ is cofinal in $D$, then $\sup D=\sup(D\cap (P_1\backslash B))=\top_1\in A$ as $A$ is Scott closed. This contradiction implies that $D\cap B$ is a cofinal subset of $D$. We now assume $D=\{(m_0,n_0,(x_k)_{a_0,b_0}): {k\in N}\}$ in the light of Remark~\ref{directed}. Note that $D\subseteq A$. Then $\sup D=(m_0,n_0,\top)\in A$. This yields that $n_0\in F$, and so $n_0\leq \max F=n^*$.
If $n_0\leq n^*-1$, then $\sup D=(m_0,n_0,\top)\in \max A\cap L_{n_0}\subseteq A\cap (\bigcup_{k\leq n^*-1}\da L_k)\subseteq A_{k_1}$. Else $n_0=n^*$. Now $D=\{(m_0,n^*,(x_k)_{a_0,b_0}): {k\in N}\}$.
For each $k\in N$,  there exists $(m_k,n_k,\top)\in (\max A\backslash (A_{k_1}\cup A_{k_2}))$ such that $(m_0,n^*, (x_k)_{a_0,b_0})\leq (m_k,n_k,\top)$. With a similar deduction as in Claim 2.1, we can get $n^*\leq n_k$. Since $n^*=\max F$ and $n_k\in F$, $n_k=n^*$.
Therefore, $(m_k,n_k,\top)=(m_k,n^*,\top)\in \max A\backslash(A_{k_1}\cup A_{k_2})$. We proceed to prove that $$\sup D=\sup_{k\in N}(m_0,n^*, (x_k)_{a_0,b_0})=(m_0,n^*, \top)\in \da(\max A\backslash (A_{k_1}\cup A_{k_2}))\cup (A\cap (\bigcup_{k\leq n^*-1}\da L_k)).$$
Again, we will only analyze the non-trivial case: when the set $\{(m_k,n^*,\top):k\in N\}$ is infinite. 
By the definition  of  $<_2$ in $P_1$, we see that $(m_0,n^*,(x_k)_{a_0,b_0})<_2 (m_k,n^*,\top)$ holds.
Observing the order of $B$ again, we know that $<_2\in \{\sqsubset_4,\ \sqsubset_1;\sqsubset_4\}$ in this case. This means that there exists $s_1\in \mathbb{N}^{<\mathbb{N}}$ satisfying that $m_0=f_{a_0,b_0}(s_1)$ and $m_k=f_{a_0,b_0}(s_1.b_k)$ for some $b_k\geq x_k$ in $\mathbb{N}$. Suppose $\sup_{k\in N}(m_0,n^*,(x_k)_{a_0,b_0})=(m_0,n^*,\top)\notin \da(\max A\backslash (A_{k_1}\cup A_{k_2}))\cup (A\cap (\bigcup_{k\leq n^*-1}\da L_k))$. Then $(m_0,n^*,\top)\in A_{k_1}\cup A_{k_2}$. Without loss of generality, we assume $(m_0,n^*,\top)\in A_{k_1}$. It indicates that $m_0=f_{a_1,b_1}(s),s\geq s_0.z_{k_1}$ from the construction of $A_{k_1}$. Since $i$ is injective, we can obtain that $(a_0,b_0)=(a_1,b_1),s=s_1$. Note that $s_1.b_k\geq s_1\geq s_0.z_{k_1}$. This reveals that $(m_k,n^*,\top)=(f_{a_0,b_0}(s_1.b_k),n^*,\top)\in A_{k_1}$, which violates the fact that $(m_k,n^*,\top)\in \max A\backslash A_{k_1}$. So indeed, Claim~$2.3$ is true. 

So
$A =  A_{k_1}\cup A_{k_2}\cup \da(\max A\backslash (A_{k_1}\cup A_{k_2}))\cup (A\cap (\bigcup_{k\leq n^*-1}\da L_k))$ is a union of three Scott closed subsets, and since $A$ is irreducible, we know that $A\subseteq A_{k_1}$, $A\subseteq A_{k_2}$ or $A\subseteq \da(\max A\backslash (A_{k_1}\cup A_{k_2}))\cup (A\cap (\bigcup_{k\leq n^*-1}\da L_k))$.

Case $1$: $A\subseteq A_{k_1}$. Then $(f_{a_1,b_1}(s_0.z_{k_2}),n^*,\top)\in A\subseteq A_{k_1}$. It follows that $(f_{a_1,b_1}(s_0.z_{k_2}),n^*,\top)\in \{(f_{a_1,b_1}(s),n^*,\top)\in \max A:s\geq s_0.z_{k_1}\}$, and so $s_0.z_{k_2}\geq s_0.z_{k_1}$. We can obtain that $s_0.z_{k_2}= s_0.z_{k_1}$ from the order on $\mathbb{N}^{<\mathbb{N}}$ and $z_{k_1},z_{k_2}\in \mathbb{N}$. Hence $z_{k_1}=z_{k_2}$. But this contradicts to our knowledge that $z_{k_1}\neq z_{k_2}$ which is stated just before~Claim 2.2.

Case $2$: $A\subseteq A_{k_2}$. This cannot be possible either, with a similar analysis as in Case $1$.

Case $3$: $A\subseteq  \da(\max A\backslash (A_{k_1}\cup A_{k_2}))\cup (A\cap (\bigcup_{k\leq n^*-1}\da L_k))$. Then $A_{k_1}\cup A_{k_2}\subseteq \da(\max A\backslash (A_{k_1}\cup A_{k_2}))\cup (A\cap (\bigcup_{k\leq n^*-1}\da L_k))$, which is absurd.

Now we have proved Claim~2. That is, $F$ is indeed infinite. 

For every $n\in F$, pick $(m_n,n,\top)\in \max A\cap L_n$. We define a function $f\in \mathbb{N}^{\mathbb{N}}$ as follows:
\begin{center}
	$$f(x)=\left\{
	\begin{aligned}
		m_x &,&x\in F\\
		1 &, & x\notin F \\
	\end{aligned}
	\right.
	$$
\end{center}

Note that $(f,n)<_3(f(n),n,\top)=(m_n,n,\top)\in A$ for every $n\in F$. Then we can conclude that $(f,n)\in A$ since $A$ is a lower set. This leads to that $\sup_{n\in F}(f,n)=\top_1\in A$ because  $A$ is Scott closed. Therefore, $A=\da \top_1$, which contradicts our assumption on that $A$ has more than one maximal elements. So our assumption that $P_1$ is not sober has been wrong, and the proof is complete. 
\end{proof}
\end{lemma}

\begin{lemma}\label{P2}
$\Sigma P_2$ is sober.
\begin{proof}
Assume, on the contrary that $\Sigma P_2$ is not sober. Then there exists an irreducible closed set $A$  of $P_2$ such that $A\notin \{\da x:x\in P_2\}$. Again,  $A=\da \max A$. For any non-trivial directed subset $D$ of $A$, since $\top_2\notin A$, $D\cap (P_2\backslash B)$ cannot be cofinal in $D$. It follows that $D\cap B$ is cofinal in $D$. As shown in the proof of Claim 1 in Lemma \ref{P1}, we also have that for any $x\in \max A$, there exists a non-trivial directed subset $D\subseteq B$, such that $x=\sup D$ and $x\in \max B$. Thus, $\max A\subseteq \max B$.

We write $F=\{n\in \mathbb{N}:\max A\cap L_n\neq \emptyset\}$, where $L_n=\mathbb{N}\times \{n\}\times L$.
Note that $D\cap B$ must be cofinal in $D$ for any non-trivial directed subset $D$ of $A$.
With a similar argument as in the proof of Claim 2 in Lemma \ref{P1}, we know that $F$ is infinite.

\textbf{Claim $1$}: There exists $n_1\geq 2$ such that $\max A\cap L_{n_1}$ is infinite.

Choose any  $(m,n,\top)\in \max A\cap L_n$, with $n\in F$ and $n\geq 2$. By Claim $1$ of Lemma \ref{P1}, there is a non-trivial directed subset $D$ of $B\cap \da(\max A\backslash \{(m,n,\top)\})$ such that $\sup D=(m,n,\top)$. We must have $D=(m,n,(x_{k})_{a,b})_{k\in N}$ for some fixed $m,n,a, b\in \mathbb{N}$, where $N$ is an infinite subset of $\mathbb{N}$. For each $k\in N$,  there exists $(m_{k},n_{k},\top)\in \max A\backslash\{(m,n,\top)\}$ such that $(m,n,(x_k)_{a,b})\leq (m_{k}, n_{k},\top)$. We note that the set $\{(m_k, n_k, \top): k\in N\}$ cannot be finite; otherwise, $\da \{(m_k, n_k, \top): k\in N\}$ is Scott closed. This implies that $\sup D=(m,n,\top)\in \da \{(m_k, n_k, \top): k\in N\}$, which means that there is $k\in N$ such that $(m,n,\top)\leq (m_k, n_k, \top)$. By the order of $P_2$, we have $(m,n,\top)= (m_k, n_k, \top)$, contradicting to the fact that $(m_{k},n_{k},\top)\in \max A\backslash\{(m,n,\top)\}$. So the set $\{(m_k, n_k, \top): k\in N\}$ is indeed infinite. Now, employing a  similar proof to that of Claim $2.1$ of Lemma \ref{P1}, we have that
$$N=\{k\in N: n_k=n\}\cup \{k\in N: n_k=n+1\}.$$
As $N$ is infinite, either $ \{k\in N: n_k=n\}$ or $\{k\in N: n_k=n+1\}$ is infinite.
As a consequence, either the set $\{(m_k, n_k, \top): n_k=n\}$ or $\{(m_k, n_k, \top): n_k=n+1\}$ is infinite.
Hence either $\max A\cap L_{n}$ or $\max A\cap L_{n+1}$ is infinite.

For every $n\in \mathbb{N}$, we set
\begin{center}
$E_n=\{m\in \mathbb{N}: \max A\cap (E_m^{a,b}\times \{n\}\times L)\neq \emptyset\}$
\end{center}

\textbf{Claim $2$}: There exists $n_2\in \mathbb{N}$ such that $E_{n_2}$ is infinite.

We prove this claim by contradiction. Assume that for any $n\in \mathbb{N}$, $E_n$ is finite, and write
\begin{itemize}
\item[] $F_1=\da (\max A\cap L_{n_1})\cup \da\{(a,n_1-1,\top),(b,n_1-1,\top):\exists ~m\in E_{n_1}, s.t.~i(a,b)=E_{m}^{a,b}\}$,
\item[] $F_2=\da (\max A\backslash L_{n_1})\cup \da\{(a,n_1,\top),(b,n_1,\top):\exists~m\in E_{n_1+1}, s.t.~i(a,b)=E_{m}^{a,b}\}$.
\end{itemize}

\textbf{Claim $2.1$}: $F_1$ is Scott closed.

It is immediate that $F_1$ is a lower set. By assumption, $E_{n_1}$ is finite. This implies that the set $$\da\{(a,n_1-1,\top),(b,n_1-1,\top):\exists ~m\in E_{n_1}, s.t.~i(a,b)=E_{m}^{a,b}\}$$ is a Scott closed subset of~$P_2$. Assume that $D$ is a non-trivial directed subset of $F_1$. 
If $D\cap \da\{(a,n_1-1,\top),(b,n_1-1,\top):\exists ~m\in E_{n_1}, s.t.~i(a,b)=E_{m}^{a,b}\}$ is cofinal in $D$, then that $\sup D \in F_1$  is trivial. Hence, we will only consider the non-trivial cases: $D\cap B \cap \da (\max A\cap L_{n_1})$ is a cofinal subset of $D$. This is because that $\top_2\notin F_1$, and hence $D\cap (P_2\setminus B)\cap \da {\max A\cap L_{n_1}}$ cannot be cofinal in $D$.
Assume $D\cap B \cap \da (\max A\cap L_{n_1})=(m_0,n_0,(x_k)_{a_0,b_0})_{k\in N}$, where $N$ is an infinite subset of $\mathbb{N}$. It follows that there exists $(m_k,n_1,\top)\in \max A\cap L_{n_1}$ such that $(m_0,n_0,(x_k)_{a_0,b_0})\leq (m_k,n_1,\top)$ for any $k\in N$.

If $n_0=n_1$, then it is easy to see that $\sup D=(m_0,n_1,\top)\in \max A\cap L_{n_1}$. Else, $n_0\neq n_1$. We will only discuss the non-trivial case: $\{(m_k,n_1,\top):k\in N\}$ is infinite. Note that both  $(m_0,n_0,(x_k)_{a_0,b_0})$ and $(m_k,n_1,\top)$ are in $B$, it follows that $(m_0,n_0,(x_k)_{a_0,b_0})<_2 (m_k,n_1,\top)$. By the definition of $ <_2$, we know that $<_2~\subseteq~ \sqsubset_2\cup\sqsubset_3\cup\sqsubset_1;\sqsubset_2\cup\sqsubset_1;\sqsubset_3$. We set
\begin{center}
	$A_{\sqsubset_{r}} = \{k\in N:(m_0,n_0,(x_k)_{a_0,b_0})\sqsubset_r (m_k,n_1,\top)\}$,
\end{center}
where $\sqsubset_{r}\in \{\ \sqsubset_{2},\ \sqsubset_{3},\ \sqsubset_{1};\sqsubset_{2},\ \sqsubset_{1};\sqsubset_{3}\}$. It follows that at least one of the four sets is infinite from the infiniteness of $N$.

Case $1$: $A_{\sqsubset_2}$ or $A_{\sqsubset_{1};\sqsubset_{2}}$ is infinite. Then $n_0=n_1-1$ and, for any $k\in A_{\sqsubset_2}$ or $A_{\sqsubset_{1};\sqsubset_{2}}$, $m_k\in E_m^{a,b}$ for some $m\in E_{n_1}$. That $E_{n_1}$ is finite guarantees the existence of $n^{*}\in E_{n_1}$ such that either $R=\{k\in A_{\sqsubset_2}:m_k\in E_{n^*}^{a,b}\}$ or $R'=\{k\in A_{\sqsubset_{1};\sqsubset_{2}}:m_k\in E_{n^*}^{a,b}\}$ is infinite. It follows that for any $k\in R\cup R'\subseteq A_{\sqsubset_2}\cup A_{\sqsubset_{1};\sqsubset_{2}}$, $m_k=f_{a_0,b_0}(y_k)$ for some $y_k\geq x_k$ in $\mathbb{N}^{<\mathbb{N}}$ and $m_0=a_0$. As $i$ is injective, we know that $(a,b)=(a_0,b_0)$. This leads to $\sup D=(m_0,n_0,\top)=(a_0,n_1-1,\top)=(a,n_1-1,\top)\in F_1$.

Case $2$: $A_{\sqsubset_3}$ or $A_{\sqsubset_{1};\sqsubset_{3}}$ is infinite. Then $n_0=n_1-1$ and, for any $k\in A_{\sqsubset_3}$ or $A_{\sqsubset_{1};\sqsubset_{3}}$, $m_k\in E_m^{a,b}$ for some $m\in E_{n_1}$. Similarly, that $E_{n_1}$ is finite guarantees the existence of $n^{*}\in E_{n_1}$ such that either $R=\{k\in A_{\sqsubset_3}:m_k\in E_{n^*}^{a,b}\}$ or $R'=\{k\in  A_{\sqsubset_{1};\sqsubset_{3}}:m_k\in E_{n^*}^{a,b}\}$ is infinite. It follows that for any $k\in R\cup R'\subseteq A_{\sqsubset_3}\cup A_{\sqsubset_{1};\sqsubset_{3}}$, $m_k=f_{a_0,b_0}(y_k)$ for some $y_k\geq x_k$ in $\mathbb{N}$ and $m_0=b_0$. As $i$ is injective, we have that $(a,b)=(a_0,b_0)$. This leads to $\sup D=(m_0,n_0,\top)=(b_0,n_1-1,\top)=(b,n_1-1,\top)\in F_1$.

Hence, $F_1$ is Scott closed.

\textbf{Claim $2.2$}: $F_2$ is Scott closed.

That the set $\da\{(a,n_1,\top),(b,n_1,\top):\exists~m\in E_{n_1+1}, s.t.~i(a,b)=E_{m}^{a,b}\}$ is a Scott closed subset of $P_2$ follows directly from the assumption that $E_{n_1+1}$ is finite. 
Let $D$ be a directed subset of $F_2$ and we proceed to show that $\sup D\in F_2$. Again, 
we will only consider the non-trivial cases: $D\cap B \cap \da (\max A\backslash  L_{n_1})$ is a cofinal subset of $D$. Assume $D\cap B \cap \da (\max A\backslash L_{n_1})=(m_0,n_0,(x_k)_{a_0,b_0})_{k\in N}$, where $N$ is an infinite subset of $\mathbb{N}$. It follows that there exists $(m_k,n_k,\top)\in \max A\backslash L_{n_1}$ such that $(m_0,n_0,(x_k)_{a_0,b_0})\leq (m_k,n_k,\top)$ for any $k\in N$.

{If $n_0\neq n_1$, then $\sup D=(m_0,n_0,\top)\in \max A\backslash L_{n_1}$. Else, $n_0=n_1$. By the order of $P_2$, we find $(m_0,n_0,(x_k)_{a_0,b_0})<_2 (m_k,n_k,\top)$. Inspecting the definition of $ <_2$ and $n_0=n_1\neq n_k$}, we have that $<_2~\subseteq~\sqsubset_2\cup\sqsubset_3\cup\sqsubset_1;\sqsubset_2\cup\sqsubset_1;\sqsubset_3$ in this case. We set
\begin{center}
	$A_{\sqsubset_{r}} = \{k\in N:(m_0,n_0,(x_k)_{a_0,b_0})\sqsubset_r (m_k,n_k,\top)\}$,
\end{center}
where $\sqsubset_{r}\in \{\ \sqsubset_{2},\ \sqsubset_{3},\ \sqsubset_{1};\sqsubset_{2},\ \sqsubset_{1};\sqsubset_{3}\}$. It follows that at least one of the four sets is infinite since $N$ is infinite. By a similar discussion as in Claim $2.1$, we can confirm that $\sup D\in F_2$ and hence that $F_2$ is Scott closed.

Note that $A\subseteq F_1\cup F_2$. Then the irreducibility of $A$ implies that $A\subseteq F_1$ or $A\subseteq F_2$. If $A\subseteq F_1$, 
then $\max A \subseteq\max F_1 \subseteq L_{n_1} \cup L_{n_1 -1}$. Note that $F$ is infinite, and this allows us to pick $n^* \in F$ with $n^* \geq n_1 +1$. It follows that $\max A \cap L_{n^*} \not= \emptyset$ by the construction of $F$. Now choose $(m^*, n^*, \top) \in \max A \cap L_{n^*}$. This means that  $(m^*, n^*, \top) \in \max A \subseteq L_{n_1} \cup L_{n_1 -1}$, a contradiction.  Else, $A\subseteq F_2$, which yields that $\max A \cap L_{n_1}\subseteq F_2$. So we  have  $$\max A \cap L_{n_1}\subseteq \{(a,n_1,\top),(b,n_1,\top):\exists~m\in E_{n_1+1}, s.t.~i(a,b)=E_{m}^{a,b}\}.$$ Note that $E_{n_1+1}$ is finite by assumption, then $\max A \cap L_{n_1}$ is finite.
This violates the fact that $\max A\cap L_{n_1}$ is infinite.

Hence, there is $n_2\in \mathbb{N}$ such that $E_{n_2}$ is infinite, and indeed \textbf{Claim $2$} is true.

For each $m\in E_{n_2}$, we pick $k_m\in E_m^{a_m,b_m}$ with $(k_m,n_2,\top)\in \max A$. For each $m\in \mathbb{N}\backslash E_{n_2}$, we pick $p_m\in E_{m}^{c_m,d_m}$. Then we define the function $f:\mathbb{N}\rightarrow \bigcup_{n\in \mathbb{N}}E_{n}^{a,b}$ as follows:
\begin{center}
	$$f(x)=\left\{
	\begin{aligned}
		k_x &,&x\in E_{n_2}\\
		p_x&, & x\in \mathbb{N}\backslash E_{n_2} \\
	\end{aligned}
	\right.
	$$
\end{center}
Note that $(f,m,n_2)<_3(f(m),n_2,\top)=(k_m,n_2,\top)\in A$ for every $m\in E_{n_2}$. As $A$ is a lower set, the set $\{(f,m,n_2):m\in E_{n_2}\}$ is a directed subset of $A$. This yields that $\sup_{m\in E_{n_2}}(f,m,n_2)=\top_2\in A$. So $A=\da \top_2$, which contradicts the fact that $A$ is not a principal ideal. Hence our assumption that $P_2$ is not sober must have been wrong. 
\end{proof}
\end{lemma}

The following result can be verified directly. We omit the proof.
\begin{proposition}\label{closed}
Let $L,M$ be two dcpo's. Then a lower subset $A$ of $L\times M$ is Scott closed if and only if for any directed subset $(x_i,y_0)_{i\in I}$ of $A$, $\sup_{i\in I}(x_i,y_0)\in A$ and any directed subset $(x_0,y_i)_{i\in I}$ of $A$, $\sup_{i\in I}(x_0,y_i)\in A$, where $x_0,y_0$ are fixed.
\end{proposition}

Finally, we arrive at our main result in this paper.

\begin{theorem}\label{productisnotsober}
$\Sigma (P_1\times P_2)$ is not sober.
\begin{proof}
Let $A=\da \{(a,a):a\in \max B\}\subseteq P_1\times P_2$.

\textcolor{black}{For the proof of the following claim, it may be useful to keep Figures~\ref{figure-6}
 and \ref{figure-7} in mind}. 
 
 \textbf{Claim $1$}: $A$ is irreducible in $\Sigma (P_1\times P_2)$.

To this end, let $U,V$ be two Scott open subsets of $\Sigma(P_1\times P_2)$ such that  $U\cap A\neq \emptyset$ and $V\cap A\neq \emptyset$. We show that $U\cap V\cap A\neq \emptyset$.

 Pick one $((m_1,n_1,\top),(m_1,n_1,\top))\in U\cap A$ and  $((m_2,n_2,\top),(m_2,n_2,\top))\in V\cap A$.  We consider  the following two cases for $n_1,n_2$.

{\bf Case $1$}: $n_1=n_2$. We only need to consider the non-trivial cases: $m_1\neq m_2$. Without loss of generality, let $m_1<m_2$. Note that $(m_2,n_2,\top)=\sup_{k\in \mathbb{N}}(m_2,n_2,(k)_{m_1,m_2})$. Then $((m_2,n_2,\top),(m_2,n_2,\top))=\sup_{k\in \mathbb{N}}((m_2,n_2,(k)_{m_1,m_2}),(m_2,n_2,(k)_{m_1,m_2}))$, which is in the Scott open set $V$. So there is  $k_0\in \mathbb{N}$ such that $((m_2,n_2,(k_0)_{m_1,m_2}),(m_2,n_2,(k_0)_{m_1,m_2}))\in V$.

By the definition of $\sqsubset_{3}$,  $(m_2,n_2,(k_0)_{m_1,m_2})\sqsubset_{3}(f_{m_{1}, m_{2}}(k_{0}),n_2+1,\top)$. As $V$ is an upper set, we know that $((f_{m_{1}, m_{2}}(k_{0}),n_2+1,\top),(f_{m_{1}, m_{2}}(k_{0}),n_2+1,\top))\in V$.

As the directed supremum 
$$\sup_{k\in \mathbb{N}}((f_{m_{1}, m_{2}}(k_{0}),n_2+1,(k)_{m_1,m_2}),
(f_{m_{1}, m_{2}}(k_{0}),n_2+1,(k)_{m_1,m_2}))$$
is equal to $((f_{m_{1}, m_{2}}(k_{0}),n_2+1,\top),(f_{m_{1}, m_{2}}(k_{0}),n_2+1,\top)$, which is in $V$.
The Scott openness of $V$ ensures that we can find $k_1\in \mathbb{N}$ such that $$((f_{m_{1},m_{2}}(k_{0}),n_2+1, (k_{1})_{m_1,m_2}),(f_{m_{1},m_{2}}(k_{0}), n_2+1, (k_{1})_{m_1,m_2}))\in V.$$
From the definition of $\sqsubset_{4}$, we have that $$(f_{m_{1},m_{2}}(k_{0}), n_2+1, (k_{1})_{m_1,m_2})\sqsubset_{4}(f_{m_{1}, m_{2}}(k_{0}.k_{1}), n_2+1,\top).$$ Again, as $V$ is an upper set,  $$((f_{m_{1},m_{2}}(k_{0}.k_{1}), n_2+1,\top), (f_{m_{1},m_{2}}(k_{0}.k_{1}), n_2+1,\top))\in V.$$

By induction on $\mathbb{N}$, for any $n\in \mathbb{N}$, there exists
\begin{center}
$((f_{m_{1},m_{2}}(k_{0}.k_{1}.\cdots .k_{n}),n_2+1,\top),(f_{m_{1},m_{2}}(k_{0}.k_{1}.\cdots .k_{n}),n_2+1,\top))\in V$.
\end{center}
The assumption that $n_1=n_2$ means that
\begin{center}
	$((f_{m_{1},m_{2}}(k_{0}.k_{1}.\cdots .k_{n}),n_1+1,\top),(f_{m_{1},m_{2}}(k_{0}.k_{1}.\cdots .k_{n}),n_1+1,\top))\in V$.
\end{center}
Note that $\{((m_1, n_1, (k_{0}.k_{1}.\cdots .k_{n})_{m_1,m_2}),(m_1, n_1, (k_{0}.k_{1}.\cdots .k_{n})_{m_1,m_2})): n\in\mathbb{N}\}$ is an increasing sequence in $P_1\times P_2$ and the directed supremum 
$$\sup \{((m_1, n_1, (k_{0}.k_{1}.\cdots .k_{n})_{m_1,m_2}),(m_1, n_1, (k_{0}.k_{1}.\cdots .k_{n})_{m_1,m_2})): n\in\mathbb{N}\}$$ 
equals $((m_1,n_1,\top),(m_1,n_1,\top))$, which is in $U$.
Thus, there is $n^*\in \mathbb{N}$ such that
$$((m_1, n_1, (k_{0}.k_{1}.\cdots .k_{n^*})_{m_1,m_2}),(m_1, n_1, (k_{0}.k_{1}.\cdots .k_{n^*})_{m_1,m_2}))\in U.$$
It turns out that $(m_1, n_1, (k_{0}.k_{1}.\cdots .k_{n^*})_{m_1,m_2}) \sqsubset_{2} (f_{m_{1},m_{2}}(k_{0}.k_1.\cdots.k_{n^*}),n_1+1,\top)$ from  the definition of $\sqsubset_{2}$. This indicates that
$$((f_{m_{1},m_{2}}(k_{0}.k_1.\cdots.k_{n^*}),n_1+1,\top),(f_{m_{1},m_{2}}(k_{0}.k_1.\cdots.k_{n^*}),n_1+1,\top))\in  U;$$
and therefore,  $$((f_{m_{1},m_{2}}(k_{0}.k_1.\cdots.k_{n^*}),n_1+1,\top),(f_{m_{1},m_{2}}(k_{0}.k_1.\cdots.k_{n^*}),n_1+1,\top))\in  U\cap V.$$
Hence, $U\cap V\cap A\neq \emptyset$. 

{\bf Case $2$}: $n_1\neq n_2$. Without loss of generality, we assume $n_1<n_2$.
The fact that 
$$((m_1, n_1, \top), (m_1, n_1, \top))\in U~\text{and}$$ 
$$((m_1, n_1, \top), (m_1, n_1, \top))=\sup_{k\in\mathbb{N}} ((m_1, n_1, (x_k)_{m_1, m_{1}+1}), (m_1, n_1, (x_k)_{m_1, m_{1}+1}))\text{,~where $x_k\in \mathbb{N}^{<\mathbb{N}}$,}$$
guarantees the existence of $k_0\in\mathbb{N}$ such that
$$((m_1, n_1, (x_{k_0})_{m_1, m_{1}+1}), (m_1, n_1, (x_{k_0})_{m_1, m_{1}+1}))\in U.$$
As $(m_1, n_1, (x_{k_0})_{m_1, m_{1}+1}) \sqsubset_2 (f_{m_1, m_1+1}(x_{k_0}), n_1+1, \top)$, we have that
 $$((f_{m_1, m_1+1}(x_{k_0}), n_1+1, \top), (f_{m_1, m_1+1}(x_{k_0}), n_1+1, \top))\in U.$$
 Thus, there is $a_1\in \mathbb{N}$ such that $((a_1, n_1+1, \top), (a_1, n_1+1, \top))\in U.$
 Applying the above arguments consecutively $n_2 - n_1$ times, we can deduce that there is $m'\in \mathbb{N}$ such that
 $$((m', n_2, \top), (m', n_2, \top))\in U.$$
 Now using the result proved in Case 1 (taking $m_1=m', n_1=n_2$), we also have $U\cap V\cap A\not=\emptyset.$

Thus, $A$ is an irreducible  set of $\Sigma (P_1\times P_2)$, and Claim~1 is true. 

\textbf{Claim $2$}: $A$ is Scott closed.

As $A$ is a lower set, by Proposition \ref{closed}, it suffices to show that for any directed subset $(x_i,y_0)_{i\in I}$ of $A$, $\sup_{i\in I}(x_i,y_0)\in A$,  and any directed subset $(x_0,y_i)_{i\in I}$ of $A$, $\sup_{i\in I}(x_0,y_i)\in A$, for any fixed $x_0,y_0$. 

\textbf{Claim $2.1$}: For any directed subset $(x_i,y_0)_{i\in I}$ of $A$, $\sup_{i\in I}(x_i,y_0)\in A$, where $y_0$ is fixed.

We only consider the non-trivial case where $(x_i,y_0)_{i\in I}$ is non-trivial. 
This implies that $(x_i)_{i\in I}$ is a non-trivial directed subset of $P_1$.

\textbf{Case $1$}: $(x_i)_{i\in I}\cap (P_1\backslash B)$ is a cofinal subset of $(x_i)_{i\in I}$. Then it follows that $(x_i)_{i\in I}\cap (P_1\backslash B)=\{(f_0,n)\}_{n\in N}$ for some fixed $f_0\in \mathbb{N}^{\mathbb{N}}$, where $N$ is an infinite subset of $\mathbb{N}$.

The fact that $\{(x_n,y_0)\}_{n\in N}=\{((f_0,n),y_0)\}_{n\in N}\subseteq A$ ensures the existence of $a_n\in \max B$ such that $((f_0,n),y_0)\leq (a_n,a_n)$ for each $n\in N$. Then $(f_0,n)< a_n$. By the order of $P_1$, we know that $(f_0,n)< _3 a_n$ or $(f_0,n)<_1;< _3 a_n$. We only take care of the non-trivial case: $\{a_n:n\in N\}$ is infinite. 
This induces  that $a_n=(f_{0}(k_n), k_n,\top)$ for some $k_n\geq n$ in $\mathbb{N}$.
Because $\{a_n:n\in N\}$ is infinite, $\{k_n:n\in N\}$ is infinite. Note that $y_0\leq a_n$ for all $n\in N$.
Now we consider the following two distinct cases for $y_0$.

\textbf{Case $1.1$}: $y_0\in P_2\backslash B$. Then let $y_0=(g_0,n_0,q_0)\in (X\times \mathbb{N})\times \mathbb{N}$. It follows that $y_0< _3 a_n$ or $y_0< _1;< _3 a_n$ in $P_2$. This yields that for each $n\in N$, $a_n=(g_0(m_n), q_0,\top)$ for some $m_n\geq n_0$ in $\mathbb{N}$. It is a contradiction to the infiniteness of $\{k_n:n\in N\}$ (note that $a_n=(f_{0}(k_n), k_n,\top)$).

\textbf{Case $1.2$}: $y_0\in B$. Then let $y_0=(m_0,n_0,(x_0)_{a_0,b_0})$. Therefore, $y_0<  a_n$ for each  $n\in N$. It turns out that $<~=~<_2~\subseteq~\sqsubset_1\cup \sqsubset_2\cup\sqsubset_3\cup\sqsubset_4\cup\sqsubset_1;\sqsubset_2\cup\sqsubset_1;\sqsubset_3\cup \sqsubset_1;\sqsubset_4$.
We set
\begin{center}
	$A_{\sqsubset_{r}} = \{k\in N:(m_0,n_0,(x_0)_{a_0,b_0})\sqsubset_r a_k=(m_k,n_k,\top)\}$,
\end{center}
where $\sqsubset_{r}\in \{\ \sqsubset_{1},\ \sqsubset_{2},\ \sqsubset_{3},\ \sqsubset_{4},\ \sqsubset_{1};\sqsubset_{2},\ \sqsubset_{1};\sqsubset_{3},\ \sqsubset_{1};\sqsubset_{4}\}$. Also let 
$$D_{\sqsubset_{r}} = \{a_k: k\in A_{\sqsubset_{r}}\}.$$
Then at least one of the  $D_{\sqsubset_{r}}$ is infinite because $\{a_k: k\in N\}$ is infinite.\\
Case $1.2.1$: $D_{\sqsubset_2}$ or $D_{\sqsubset_{1};\sqsubset_{2}}$ is infinite. Then $n_k=n_0+1$ for any $k\in A_{\sqsubset_{2}}$ or $k\in A_{\sqsubset_{1};\sqsubset_{2}}$. 
Note that $a_{n}=(f_0(k_n), k_n, \top)$ in Case 1. Then $a_{n}=a_{n'}$ if $k_n=k_{n'}$. 
Thus, $D_{\sqsubset_2}=\{a_{k}:k\in A_{\sqsubset_2}\}$ and $D_{\sqsubset_{1};\sqsubset_{2}}=\{a_k: k\in A_{\sqsubset_{1};\sqsubset_{2}}\}$ are finite, a contradiction.\\
Case $1.2.2$: $D_{\sqsubset_3}$ or $D_{\sqsubset_{1};\sqsubset_{3}}$ is infinite. Then $n_k=n_0+1$ for any $k\in A_{\sqsubset_{3}}$ or $k\in A_{\sqsubset_{1};\sqsubset_{3}}$. Note that $a_{n}=(f_0(k_n), k_n, \top)$ in Case 1. Then $a_{n}=a_{n'}$ if $k_n=k_{n'}$. 
Thus, $D_{\sqsubset_3}=\{a_{k}:k\in A_{\sqsubset_3}\}$ and $D_{\sqsubset_{1};\sqsubset_{3}}=\{a_k: k\in A_{\sqsubset_{1};\sqsubset_{3}}\}$ are finite, a contradiction.\\
\textcolor{black}{
Case $1.2.3$: $D_{\sqsubset_4}$ or $D_{\sqsubset_{1};\sqsubset_{4}}$ or $D_{\sqsubset_{1}}$ is infinite. Then $n_k=n_0$ for all $k\in A_{\sqsubset_{4}}$ or $k\in A_{\sqsubset_{1};\sqsubset_{4}}$ or $k\in A_{\sqsubset_{1}}$, a similar reasoning as in the above two cases leads to absurdness.}

\textbf{Case $2$}: $(x_i)_{i\in I}\cap B$ is a cofinal subset of $(x_i)_{i\in I}$. 

Then we assume $(x_i)_{i\in I}\cap B=(m_0,n_0,(x_k)_{a_0,b_0})_{k\in N}$ for fixed natural numbers $m_0,n_0,a_0$ and $b_0$, where $N$ is an infinite subset of $\mathbb{N}$. We can find $a_k\in \max B$ such that $((m_0,n_0,(x_k)_{a_0,b_0}), y_0)< (a_k,a_k)$. We will only consider the non-trivial case: when $\{a_k:k\in N\}$ is infinite. By the order of $P_1$, we know that $((m_0,n_0,(x_k)_{a_0,b_0}), y_0)<_2 (a_k,a_k)$. It follows that $<_2\subseteq\sqsubset_1\cup \sqsubset_2\cup\sqsubset_3\cup\sqsubset_4\cup\sqsubset_1;\sqsubset_2\cup\sqsubset_1;\sqsubset_3\cup \sqsubset_1;\sqsubset_4$ in this case. We set
\begin{center}
	$A_{\sqsubset_{r}} = \{k\in N:(m_0,n_0,(x_k)_{a_0,b_0})\sqsubset_r a_k=(m_k,n_k,\top)\}$,
\end{center}
where $\sqsubset_{r}\in \{\ \sqsubset_{1},\ \sqsubset_{2},\ \sqsubset_{3},\ \sqsubset_{4},\ \sqsubset_{1};\sqsubset_{2},\ \sqsubset_{1};\sqsubset_{3},\ \sqsubset_{1};\sqsubset_{4}\}$. 
\textcolor{black}{Similarly, we define $D_{\sqsubset_{r}} =  \{a_k: k\in A_{\sqsubset_{r}}\}$.} Then at least one of the $D_{\sqsubset_{r}}$ is infinite. If $D_{\sqsubset_{1}}$ is infinite, then $a_n=(m_0,n_0,\top)$ for each $a_n\in D_{\sqsubset_{1}}$, which is fixed. This contradicts to the assumption that $D_{\sqsubset_{1}}$ is infinite. Now we distinguish the following cases:

\textbf{Case $2.1$}: $D_{\sqsubset_2}$ or $D_{\sqsubset_{1};\sqsubset_{2}}$ is infinite. Then for any $k\in A_{\sqsubset_2}$ or $k\in A_{\sqsubset_{1};\sqsubset_{2}}$, $n_k=n_0+1$, $m_k=f_{a_0,b_0}(y_k)$ for some $y_k\geq x_k$ in $\mathbb{N}^{<\mathbb{N}}$, and $m_0=a_0$. This yields that $m_k\in i(a_0,b_0)=E_{n^*}^{a_0,b_0}$ for some $n^*\in \mathbb{N}$.
Also note that $\{x_k\}_{k\in N}$ is a non-trivial directed subset of $\mathbb{N}^{<\mathbb{N}}$, hence the \textcolor{black}{lengths of $x_{k} (k\in N)$ range over an infinite subset of $\mathbb{N}$, and so do the lengths of $y_k (k\in N)$.}

Now we consider the following two distinct subcases for $y_0$.

Case $2.1.1$: $y_0\in P_2\backslash B$. Then assume $y_0=(g_0,n_0,k_0)$. By the order of $P_2$, we know that $y_0<_3 a_k$ or $y_0<_1;<_3 a_k$. This implies that $a_k=(g_0(p_k),k_0,\top)$ for some $p_k\geq n_0$ in $\mathbb{N}$. That either $D_{\sqsubset_2}$ or $D_{\sqsubset_{1};\sqsubset_{2}}$ is infinite ensures that either $\{p_k: k\in A_{\sqsubset_2}\}$ or $\{p_k:k\in A_{\sqsubset_1;\sqsubset_2}\}$ is infinite. Now we know that $m_k=g_0(p_k)\in E_{p_k}^{c,d}$ and $m_k\in i(a_0,b_0)=E_{n^*}^{a_0,b_0}$ from Case 2.1.  So $E_{p_k}^{c,d}=E_{n^*}^{a_0, b_0}$ by the item (1) of Remark~\ref{i}, and then $p_k=n^*$ for all $k\in A_{\sqsubset_2}$ or $k\in A_{\sqsubset_{1};\sqsubset_{2}}$. But it contradicts the fact that either $\{p_k: k\in A_{\sqsubset_2}\}$ or $\{p_k:k\in A_{\sqsubset_1;\sqsubset_2}\}$ is infinite.

Case $2.1.2$: $y_0\in B$. Then assume $y_0=(m^*,n^*,(x^*)_{a_*,b_*})\in P_2$. It turns out that $y_0< a_k$ for each $k\in N.$ By the order of $P_2$, we know that $<=<_2\subseteq \sqsubset_1\cup \sqsubset_2\cup\sqsubset_3\cup\sqsubset_4\cup\sqsubset_1;\sqsubset_2\cup\sqsubset_1;\sqsubset_3\cup \sqsubset_1;\sqsubset_4$.
We set
\begin{center}
	$B_{\sqsubset_{r}} = \{k\in A_{\sqsubset_2}~\mathrm{or}~k\in A_{\sqsubset_1;\sqsubset_2}:y_0\sqsubset_r a_k=(m_k,n_k,\top)\}$,
\end{center}
where $\sqsubset_{r}~\in \{\sqsubset_{1},\ \sqsubset_{2},\ \sqsubset_{3},\ \sqsubset_{4},\ \sqsubset_{1};\sqsubset_{2},\ \sqsubset_{1};\sqsubset_{3},\ \sqsubset_{1};\sqsubset_{4}\}$. It follows that at least one of the seven sets is infinite as $A_{\sqsubset_2}$ or $A_{\sqsubset_1;\sqsubset_2}$ is infinite. If $B_{\sqsubset_{1}}$ is infinite, then for each $k\in B_{\sqsubset_{1}}$, $a_k=(m^*,n^*,\top)$, which is fixed. Note that in Case 2.1, \textcolor{black}{the lengths of $x_k (k\in N)$ form an infinite set} and $y_k\geq x_k$ in $\mathbb{N}^{<\mathbb{N}}$. Thus, \textcolor{black}{the lengths of $y_k (k\in B_{\sqsubset_{1}})$  also form an infinite set.} This means that $\{a_k:k\in B_{\sqsubset_{1}} \}$ is infinite, which is a contradiction to the fact that $a_k=(m^*,n^*,\top)$ for  $k\in B_{\sqsubset_{1}}$.

Case $2.1.2.1$: $B_{\sqsubset_2}$ or $B_{\sqsubset_{1};\sqsubset_{2}}$ is infinite. Then for any $k\in B_{\sqsubset_2}$ or $k\in B_{\sqsubset_{1};\sqsubset_{2}}$, $m_k=f_{a^*,b^*}(z_k)$ for some $z_k\geq x^*$ in $\mathbb{N}^{<\mathbb{N}}$ and $m^*=a^*$, $n^*+1=n_k$, $n_k=n_0+1$ from Case 2.1, and so $n^*+1=n_k=n_0+1$. Note that $m_k=f_{a_0,b_0}(y_k)$ for any $k\in A_{\sqsubset_2}$ or $k\in A_{\sqsubset_{1};\sqsubset_{2}}$ from Case 2.1. Then we have that $(a^*,b^*)=(a_0,b_0)$ as $i$ is injective. Therefore,
\begin{align*}
    \sup_{i\in I}(x_i,y_0)& =\sup_{k\in N}((m_0,n_0,(x_k)_{a_0,b_0}),y_0)\\
    & =((m_0,n_0,\top),y_0)\\
    & =((a_0,n_0,\top),(a_0,n_0,(x^*)_{a_0,b_0}))\\
    & \in \da((a_0,n_0,\top), (a_0,n_0,\top))\\
    & \subseteq A.
\end{align*}

Case $2.1.2.2$: $B_{\sqsubset_3}$ or $B_{\sqsubset_{1};\sqsubset_{3}}$ is infinite. Then for each $k\in B_{\sqsubset_3}$ or $k\in B_{\sqsubset_{1};\sqsubset_{3}}$, $m_k=f_{a^*,b^*}(z_k)$ for some $z_k\geq x^*$ in $\mathbb{N}$ and $m^*=b^*$, $n^*+1=n_k$, $n_k=n_0+1$ from Case 2.1, and so $n^*+1=n_k=n_0+1$.
Note that $m_k=f_{a_0,b_0}(y_k)$ for any $k\in A_{\sqsubset_2}$ or $k\in A_{\sqsubset_{1};\sqsubset_{2}}$ by Case 2.1. Then due to the injectivity of $i$, we know that $(a^*,b^*)=(a_0,b_0)$, $y_k=z_k$. Note that the length of $z_k$ is 1.

Also, In Case 2.1, we get that the \textcolor{black}{lengths of $x_k (k\in N)$ form an infinite set} and $y_k\geq x_k$ in $\mathbb{N}^{<\mathbb{N}}$ for all $k\in A_{\sqsubset_2}$ or $k\in A_{\sqsubset_{1};\sqsubset_{2}}$. Hence, \textcolor{black}{the lengths of $y_k (k\in B_{\sqsubset_3})$ or that of $y_k (k\in B_{\sqsubset_{1};\sqsubset_{3}})$ form infinite sets.}
Then there is $k_0\in B_{\sqsubset_3}$ or $k_0\in B_{\sqsubset_{1};\sqsubset_{3}}$ such that $y_{k_0}$ whose length is greater than $1$.
This contradicts the results that $y_{k_0}=z_{k_0}$ and the length of $z_{k_0}$ equals to 1.

Case $2.1.2.3$: $B_{\sqsubset_4}$ or $B_{\sqsubset_{1};\sqsubset_{4}}$ is infinite. Then for each $k\in B_{\sqsubset_4}$ or $k\in B_{\sqsubset_{1};\sqsubset_{4}}$, there exists $s^*\in \mathbb{N}^{<\mathbb{N}}$ such that $m^*=f_{a^*,b^*}(s^*)$, $m_k=f_{a^*,b^*}(s^*.z_k)$ for some $z_k\geq x^*$ in $\mathbb{N}$, and $n_0+1=n_k=n^*$. Note that $m_k=f_{a_0,b_0}(y_k)$ for any $k\in A_{\sqsubset_2}$ or $k\in A_{\sqsubset_{1};\sqsubset_{2}}$ from Case 2.1. Then the injectivity of $i$ implies that $(a^*,b^*)=(a_0,b_0)$ and $y_k=s^*.z_k$. This yields that the length of string $y_k$ is fixed as $\mid s^*\mid +1$. By a similar discussion for the infiniteness of the \textcolor{black}{set of the lengths of $y_k (k\in B_{\sqsubset_3})$ or $y_k (k\in B_{\sqsubset_{1};\sqsubset_{3}})$ in Case 2.1.2.2,} we can obtain that the \textcolor{black}{lengths of $y_k (k\in B_{\sqsubset_4})$ or that of $y_k (k\in B_{\sqsubset_{1};\sqsubset_{4}})$ form infinite sets}, which is a contradiction to the fact that \textcolor{black}{the lengths of the strings $y_k (k\in B_{\sqsubset_4})$ or $y_k (k\in B_{\sqsubset_{1};\sqsubset_{4}})$ are fixed as $\mid s^*\mid +1$.}

\textbf{Case $2.2$}: $D_{\sqsubset_3}$ or $D_{\sqsubset_{1};\sqsubset_{3}}$ is infinite. Then for each  $k\in A_{\sqsubset_3}$ or $k\in A_{\sqsubset_{1};\sqsubset_{3}}$, $n_k=n_0+1$, $m_k=f_{a_0,b_0}(y_k)$ for some $y_k\geq x_k$ in $\mathbb{N}$, and $m_0=b_0$. This leads to that $m_k\in i(a_0,b_0)=E_{n^*}^{a_0,b_0}$ for some $n^*\in \mathbb{N}$. 
Again, $\{y_k\}_{k\in N}$ should be an infinite set of $\mathbb{N}$.

We consider the two distinguished cases for $y_0$.

Case $2.2.1$: $y_0\in P_2\backslash B$. We assume $y_0=(g_0,n_0,k_0)$. By the order of $P_2$, we have that $y_0<_3 a_k$ or $y_0<_1;<_3 a_k$.
This means that $a_k=(g_0(p_k),k_0,\top)$ for some $p_k\geq n_0$ in $\mathbb{N}$. The infiniteness of $D_{\sqsubset_3}$ or $D_{\sqsubset_{1};\sqsubset_{3}}$  ensures the infiniteness of $\{p_k:k\in A_{\sqsubset_3}\}$ or $\{p_k:k\in A_{\sqsubset_1;\sqsubset_3}\}$, whence, $m_k=g_0(p_k)\in E_{p_k}^{c,d}$, and $m_k\in E_{n^*}^{a_0,b_0}$ from Case 2.2. This implies that $p_k=n^*$ for all  $k\in A_{\sqsubset_3}$ or $k\in A_{\sqsubset_1;\sqsubset_3}$, a contradiction to the infiniteness of $\{p_k:k\in A_{\sqsubset_3}\}$ or $\{p_k:k\in A_{\sqsubset_1;\sqsubset_3}\}$.

Case $2.2.2$: $y_0\in B$. Then assume $y_0=(m^*,n^*,(x^*)_{a_*,b_*})\in P_2$. It follows that $y_0< a_k$ for any $k\in N.$ Observing the order of $P_2$, we have that $<=<_2\subseteq\sqsubset_1\cup \sqsubset_2\cup\sqsubset_3\cup\sqsubset_4\cup\sqsubset_1;\sqsubset_2\cup\sqsubset_1;\sqsubset_3\cup \sqsubset_1;\sqsubset_4$ in this case. 
We set
\begin{center}
	$B_{\sqsubset_{r}} = \{k\in A_{\sqsubset_3}~\mathrm{or}~k\in A_{\sqsubset_1;\sqsubset_3}:y_0\sqsubset_r a_k=(m_k,n_k,\top)\}$,
\end{center}
where $\sqsubset_{r}\in \{\ \sqsubset_{1},\sqsubset_{2},\ \sqsubset_{3},\ \sqsubset_{4},\ \sqsubset_{1};\sqsubset_{2},\ \sqsubset_{1};\sqsubset_{3},\ \sqsubset_{1};\sqsubset_{4}\}$. Again,  at least one of the seven sets is infinite by the infiniteness of $A_{\sqsubset_3}$ or $A_{\sqsubset_1;\sqsubset_3}$. If $B_{\sqsubset_{1}}$ is infinite, then for each $k\in B_{\sqsubset_{1}}$, $a_k=(m^*,n^*,\top)$. Note that in Case 2.2, $\{x_k :k\in N\}$ is an infinite set of $\mathbb{N}$ and $y_k\geq x_k$ in $\mathbb{N}$. Thus $\{y_k :k\in B_{\sqsubset_{1}}\}$ is also an infinite set of $\mathbb{N}$. This means that $\{a_k:k\in B_{\sqsubset_{1}} \}$ is infinite, which is a contradiction to the fact that $a_k=(m^*,n^*,\top)$ for any $k\in B_{\sqsubset_{1}}$.

Case $2.2.2.1$: $B_{\sqsubset_2}$ or $B_{\sqsubset_{1};\sqsubset_{2}}$ is infinite. Then for each $k\in B_{\sqsubset_2}$ or $k\in B_{\sqsubset_{1};\sqsubset_{2}}$, $m_k=f_{a^*,b^*}(z_k)$ for some $z_k\geq x^*$ in $\mathbb{N}^{<\mathbb{N}}$ and $m^*=a^*$, $n^*+1=n_k=n_0+1$. Note that $m_k=f_{a_0,b_0}(y_k)$ for any $k\in A_{\sqsubset_3}$ or $k\in A_{\sqsubset_{1};\sqsubset_{3}}$ from Case 2.2. As $i$ is injective, we know that $(a^*,b^*)=(a_0,b_0)$ and $z_k=y_k\geq x^*$ in $\mathbb{N}^{<\mathbb{N}}$. But the length of $y_k$ is 1, it holds that $y_k=x^*$, which is fixed. Note that in Case 2.2, $\{x_k :k\in N\}$ is an infinite set of $\mathbb{N}$ and $y_k\geq x_k$ in $\mathbb{N}$. Thus, $\{y_k :k\in B_{\sqsubset_{2}}\}$ or \{$y_k: k\in B_{\sqsubset_{1};\sqsubset_{2}}\}$ is also an infinite set of $\mathbb{N}$. This contradicts the fact that $y_k=x^*$ for any $k\in B_{\sqsubset_{2}}$ or $k\in B_{\sqsubset_{1};\sqsubset_{2}}$.

Case $2.2.2.2$: $B_{\sqsubset_3}$ or $B_{\sqsubset_{1};\sqsubset_{3}}$ is infinite. Then for each $k\in B_{\sqsubset_3}$ or $k\in B_{\sqsubset_{1};\sqsubset_{3}}$, $m_k=f_{a^*,b^*}(z_k)$ for some $z_k\geq x^*$ in $\mathbb{N}$ and $m^*=b^*$, $n^*+1=n_k$, $n_k=n_0+1$ by Case 2.2, and so $n^*+1=n_k=n_0+1$.
Note that $m_k=f_{a_0,b_0}(y_k)$ for each $k\in A_{\sqsubset_3}$ or $k\in A_{\sqsubset_{1};\sqsubset_{3}}$ from Case 2.2. Then due to the injectivity of $i$, we know that $(a^*,b^*)=(a_0,b_0)$, $y_k=z_k$. It follows that $m^*=b^*=b_0=m_0$ and $n^*=n_0$. This implies that $\sup_{i\in I}(x_i,y_0)=((m_0,n_0,\top), (m^*,n^*,(x^*)_{a^*,b^*}))=((m_0,n_0,\top), (m_0,n_0,(x^*)_{a^*,b^*}))\in \da ((m_0,n_0,\top), (m_0,n_0,\top))\subseteq A$.

Case $2.2.2.3$: $B_{\sqsubset_4}$ or $B_{\sqsubset_{1};\sqsubset_{4}}$ is infinite. Then for each $k\in B_{\sqsubset_4}$ or $k\in B_{\sqsubset_{1};\sqsubset_{4}}$, there exists $s^*\in \mathbb{N}^{<\mathbb{N}}$ such that $m^*=f_{a^*,b^*}(s^*)$, $m_k=f_{a^*,b^*}(s^*.z_k)$ for some $z_k\geq x^*$ in $\mathbb{N}$ and $n_0+1=n_k=n^*$. Note that $m_k=f_{a_0,b_0}(y_k)$ for any $k\in A_{\sqsubset_3}$ or $k\in A_{\sqsubset_{1};\sqsubset_{3}}$ from Case 2.2. Then the injectivity of $i$ implies that $(a^*,b^*)=(a_0,b_0)$ and $y_k=s^*.z_{k}$. This yields that the length of string $y_k$ is fixed as $\mid s^*\mid +1\geq 2$, which contradicts that $y_k\in \mathbb{N}$ from Case 2.2.

\textbf{Case $2.3$}: $D_{\sqsubset_4}$ or $D_{\sqsubset_{1};\sqsubset_{4}}$ is infinite. Then for each $k\in A_{\sqsubset_4}$ or $k\in A_{\sqsubset_{1};\sqsubset_{4}}$, there exists $s_0\in \mathbb{N}^{<\mathbb{N}}$ such that $m_0=f_{a_0,b_0}(s_0)$, $m_k=f_{a_0,b_0}(s_0.y_k)$ for some $y_k\geq x_k$ in $\mathbb{N}$ and $n_k=n_0$. This yields that $m_k\in i(a_0,b_0)=E_{n^*}^{a_0,b_0}$ for some $n^*\in \mathbb{N}$. Again, the set $\{y_k:k\in N\}$ should be an infinite set of $\mathbb{N}$ since $\{x_k:k\in N\}$ is an infinite subset of $\mathbb{N}$.
Now we consider two distinct cases for $y_0$.

Case $2.3.1$: $y_0\in P_2\backslash B$. Assume $y_0=(g_0,n_0,k_0)$. By the order of $P_2$, it follows that $y_0<_3 a_k$ or $y_0<_1;<_3 a_k$.
This means that $a_k=(g_0(p_k),k_0,\top)$ for some $p_k\geq n_0$ in $\mathbb{N}$. The infiniteness of $D_{\sqsubset_4}$ or $D_{\sqsubset_{1};\sqsubset_{4}}$ ensures that either  $\{p_k:k\in A_{\sqsubset_4}\}$ or $\{p_k:k\in A_{\sqsubset_1;\sqsubset_4}\}$ is infinite, whence $m_k=g_0(p_k)\in E_{p_k}^{c,d}$. Note that $m_k\in E_{n^*}^{a_0,b_0}$ by Case 2.3, then $m_k\in E_{p_k}^{c,d}\cap E_{n^*}^{a_0,b_0}$.
This leads to that $p_k=n^*$ for any $k\in A_{\sqsubset_4}$ or $k\in A_{\sqsubset_1;\sqsubset_4}$, a contradiction to the infiniteness of $\{p_k:k\in A_{\sqsubset_4}\}$ or $\{p_k:k\in A_{\sqsubset_1;\sqsubset_4}\}$.

Case $2.3.2$: $y_0\in B$. Then assume $y_0=(m^*,n^*,(x^*)_{a^*,b^*})\in P_2$. It turns out that $y_0<_2 a_k$ for any $k\in N.$ By the order of $P_2$, we have that $<_2\subseteq\sqsubset_1\cup \sqsubset_2\cup\sqsubset_3\cup\sqsubset_4\cup\sqsubset_1;\sqsubset_2\cup\sqsubset_1;\sqsubset_3\cup \sqsubset_1;\sqsubset_4$. We set
\begin{center}
	$B_{\sqsubset_{r}} = \{k\in A_{\sqsubset_4}~\mathrm{or}~k\in A_{\sqsubset_1;\sqsubset_4}:y_0\sqsubset_r a_k=(m_k,n_k,\top)\}$,
\end{center}
where $\sqsubset_{r}\in \{\ \sqsubset_{1},\ \sqsubset_{2},\ \sqsubset_{3},\ \sqsubset_{4},\ \sqsubset_{1};\sqsubset_{2},\ \sqsubset_{1};\sqsubset_{3},\ \sqsubset_{1};\sqsubset_{4}\}$. Again, at least one of the seven sets is infinite by the infiniteness of $A_{\sqsubset_4}$ or $A_{\sqsubset_1;\sqsubset_4}$. If $B_{\sqsubset_{1}}$ is infinite, then for any $k\in B_{\sqsubset_{1}}$, $a_k=(m^*,n^*,\top)$, which is fixed. Note that in Case 2.3, the set $\{x_k : k\in N\}$ is an infinite set of $\mathbb{N}$ and $y_k\geq x_k$ in $\mathbb{N}$. Thus, $\{y_k :k\in B_{\sqsubset_{1}}\}$ is also an infinite set of $\mathbb{N}$. This means that $\{a_k:k\in B_{\sqsubset_{1}} \}$ is infinite, which is a contradiction to the fact that $a_k=(m^*,n^*,\top)$ for any $k\in B_{\sqsubset_{1}}$.

Case $2.3.2.1$: $B_{\sqsubset_2}$ or $B_{\sqsubset_{1};\sqsubset_{2}}$ is infinite. Then for each $k\in B_{\sqsubset_2}$ or $k\in B_{\sqsubset_{1};\sqsubset_{2}}$, $m_k=f_{a^*,b^*}(z_k)$ for some $z_k\geq x^*$ in $\mathbb{N}^{<\mathbb{N}}$ and $m^*=a^*$, $n^*+1=n_k,$ $n_k=n_0$ from Case 2.3, and so $n^*+1=n_k=n_0$. By using similar argument as above, we can obtain that $\{y_k :k\in B_{\sqsubset_{2}}\}$ or $\{y_k :k\in B_{\sqsubset_{1};\sqsubset_{2}}\}$ is also an infinite subset of $\mathbb{N}$.  
Note that $m_k=f_{a_0,b_0}(s_0.y_k)$ for any $k\in A_{\sqsubset_4}$ or $k\in A_{\sqsubset_{1};\sqsubset_{4}}$ by Case 2.3. Then by the injectivity of $i$, we know that $(a^*,b^*)=(a_0,b_0)$ and $z_k=s_0.y_k\geq x^*$ for any $k\in B_{\sqsubset_2}$ or $k\in B_{\sqsubset_{1};\sqsubset_{2}}$. We obtain that $x^*\leq s_0$ in $\mathbb{N}^{<\mathbb{N}}$ as $\{y_k:k\in B_{\sqsubset_2}\}$ or $\{y_k:k\in B_{\sqsubset_{1};\sqsubset_{2}}\}$ is infinite, and $y_k\in \mathbb{N}$. It follows that $\sup_{i\in I}(x_i,y_0)=((m_0,n_0,\top),(m^*,n^*,(x^*)_{a^*,b^*}))=((f_{a_0,b_0}(s_0),n_0,\top),(a_0,n_0-1,(x^*)_{a_0,b_0}))$. Since $(a_0,n_0-1,(x^*)_{a_0,b_0}))\sqsubset_1 (a_0,n_0-1,(s_0)_{a_0,b_0}))\sqsubset_2(f_{a_0,b_0}(s_0),n_0,\top)$, we know that $\sup_{i\in I}(x_i,y_0)\in \da ((f_{a_0,b_0}(s_0),n_0,\top),(f_{a_0,b_0}(s_0),n_0,\top))\subseteq  A$.

Case $2.3.2.2$: $B_{\sqsubset_3}$ or $B_{\sqsubset_{1};\sqsubset_{3}}$ is infinite. Then for each $k\in B_{\sqsubset_3}$ or $k\in B_{\sqsubset_{1};\sqsubset_{3}}$, $m_k=f_{a^*,b^*}(z_k)$ for some $z_k\geq x^*$ in $\mathbb{N}$ and $m^*=b^*$, $n^*+1=n_k=n_0$.
Note that $m_k=f_{a_0,b_0}(s_0.y_k)$ for any $k\in A_{\sqsubset_4}$ or $k\in A_{\sqsubset_{1};\sqsubset_{4}}$ by Case 2.3. Then by the injectivity of $i$,  $(a^*,b^*)=(a_0,b_0)$, $s_0.y_k=z_k$, which is a contradiction to the fact that $z_k$ is a natural number.

Case $2.3.2.3$: $B_{\sqsubset_4}$ or $B_{\sqsubset_{1};\sqsubset_{4}}$ is infinite. Then for each $k\in B_{\sqsubset_4}$ or $k\in B_{\sqsubset_{1};\sqsubset_{4}}$, there exists $s^*\in \mathbb{N}^{<\mathbb{N}}$ such that $m^*=f_{a^*,b^*}(s^*)$, $m_k=f_{a^*,b^*}(s^*.z_k)$ for some $z_k\geq x^*$ in $\mathbb{N}$, $n_k=n^*$, and $n_0=n_k$ from Case 2.3. So $n_0=n_k=n^*$. Note that $m_k=f_{a_0,b_0}(s_0.y_k)$ for any $k\in A_{\sqsubset_4}$ or $k\in A_{\sqsubset_{1};\sqsubset_{4}}$ by Case 2.3. Then the injectivity of $i$ implies that $(a^*,b^*)=(a_0,b_0)$ and $s_0.y_k=s^*.z_k$. It follows that $s_0=s^*$ because $y_k, z_k\in\mathbb{N}$. Hence, $m_0=f_{a_0,b_0}(s_0)=f_{a^*,b^*}(s^*)=m^*$. This yields that 
\begin{align*}
    \sup_{i\in I}(x_i,y_0) & =((m_0,n_0,\top), (m^*,n^*,(x^*)_{a^*,b^*}))\\ 
    & =((m_0,n_0,\top),(m_0,n_0,(x^*)_{a^*,b^*}))\\
    & \in  \da ((m_0,n_0,\top),(m_0,n_0,\top))\\
    & \subseteq A.
\end{align*}

All the above cases show that Claim 2.1 is true. 

\textbf{Claim $2.2$}:  For every directed subset $(x_0,y_i)_{i\in I}$ of $A$, $\sup_{i\in I}(x_0,y_i)\in A$, where $x_0$ is fixed.
Then $(y_i)_{i\in I}$ is a directed subset of $P_2$. We only consider the case when $(y_i)_{i\in I}$ is non-trivial.

\textbf{Case $1$}: $(y_i)_{i\in I}\cap (P_2\backslash B)$ is cofinal. Then we assume that $(y_i)_{i\in I}\cap (P_2\backslash B)=(g_0,n,k_0)_{n\in N}$, where $N$ is an infinite subset of $\mathbb{N}$. It follows that there exists $a_n\in \max B$ such that $(x_0,(g_0,n,k_0))\leq (a_n,a_n)$. We need to prove that $\sup_{i\in I}(x_0,y_i)=\sup_{n\in N}(x_0,(g_0,n,k_0))\in A$.
We will only discuss the non-trivial case: $\{a_n:n\in N\}$ is infinite.
Then $(g_0,n,k_0)<a_n$ for all $n\in N$. By the order of $P_2$, we know that $(g_0,n,k_0)<_1;<_3 a_n$ or $(g_0,n,k_0)<_3 a_n$ holds. This means that $a_n=(g_0(p_n),k_0,\top)$ for some $p_n\geq n$. So we have $\{p_n:n\in N\}$ is infinite. Now we consider the following two subcases for $x_0$.

\textbf{Case $1.1$}: $x_0\in P_1\backslash B$. Then assume $x_0=(f_0,n_0)$. Note that $x_0\leq a_n$ for all $n\in N$. By the order of $P_1$, we have that $x_0<_3 a_n$ or $x_0<_1;<_3 a_n$ in $P_1$. This means that $a_n=(f_0(k_n),k_n,\top)$ for some $k_n\geq n_0$. The infiniteness of $\{a_n:n\in N\}$ implies that $\{k_n:n\in N\}$ is infinite, which is a contradiction to the fact  that $a_n=(g_0(p_n),k_0,\top)$ for all $n\in N$.

\textbf{Case $1.2$}: $x_0\in B$. Assume that  $x_0=(m^*,n^*,(x^*)_{a^*,b^*})$. It turns out that $x_0< a_n$ for all $n\in N$. Following the definition of $<$, we know that $<=<_2\subseteq\sqsubset_1\cup \sqsubset_2\cup\sqsubset_3\cup\sqsubset_4\cup\sqsubset_1;\sqsubset_2\cup\sqsubset_1;\sqsubset_3\cup\sqsubset_1;\sqsubset_4$. 
Again, let
\begin{center}
	$A_{\sqsubset_{r}} = \{n\in N:(m^*,n^*,(x^*)_{a^*,b^*})\sqsubset_r a_n\}$,
\end{center}
where $\sqsubset_{r}\in \{\ \sqsubset_{1},\ \sqsubset_{2},\ \sqsubset_{3},\ \sqsubset_4,\ \sqsubset_{1};\sqsubset_{2},\ \sqsubset_{1};\sqsubset_{3},\ \sqsubset_1;\sqsubset_4\}$. 
Also let $$D_{\sqsubset_{r}} =\{a_k: k\in A_{\sqsubset_{r}}\}.$$

At least one of the $D_{\sqsubset_{r}}$ is infinite by the infiniteness of $\{a_n:n\in N\}$.

Case $1.2.1$: $D_{\sqsubset_2}$ or $D_{\sqsubset_{1};\sqsubset_{2}}$ is infinite. Then $a_n=(f_{a^*,b^*}(y_n),n^*+1,\top)$ for some $y_n\geq x^*$ in $\mathbb{N}^{<\mathbb{N}}$, for any $k\in A_{\sqsubset_2}$ or $k\in A_{\sqsubset_{1};\sqsubset_{2}}$. It follows that $i(a^*,b^*)=E_{r_0}^{a^*,b^*}$ for some fixed $r_0\in \mathbb{N}$.
This means that $\{a_n:n\in A_{\sqsubset_2}\}\subseteq E_{r_0}^{a^*,b^*}\times \{n^*+1\}\times L$ or $\{a_n:n\in A_{\sqsubset_1;\sqsubset_2}\}\subseteq E_{r_0}^{a^*,b^*}\times \{n^*+1\}\times L$ for some fixed $r_0\in \mathbb{N}$. Note that $a_n=(g_0(p_n),k_0,\top)$ for any $n\in N$ by Case 1. Then $g_0(p_n)\in E_{p_n}^{c,d}$ for all $n\in N$, which implies  $p_n=r_0$ for all $n\in N$, then $a_n$ is actually fixed as $(g_0(r_0), k_0,\top)$.  This contradicts the assumption that $D_{\sqsubset_2}$ or $D_{\sqsubset_{1};\sqsubset_{2}}$ is infinite.

Case $1.2.2$: $D_{\sqsubset_3}$ or $D_{\sqsubset_{1};\sqsubset_{3}}$ is infinite. Then $a_n=(f_{a^*,b^*}(y_n),n^*+1,\top)$ for some $y_n\geq x^*$ in $\mathbb{N}$, for each  $n\in A_{\sqsubset_3}$ or $n\in A_{\sqsubset_{1};\sqsubset_{3}}$. 

As in Case 1.2.1, we then have that $p_n=r_0$ is a constant, and so $a_n$ is fixed, contradicting  the assumption that $D_{\sqsubset_3}$ or $D_{\sqsubset_{1};\sqsubset_{3}}$ is infinite.

Case $1.2.3$: $D_{\sqsubset_4}$ or $D_{\sqsubset_{1};\sqsubset_{4}}$ is infinite. Then there exists $s^*\in \mathbb{N}^{<\mathbb{N}}$ such that $m^*=f_{a^*,b^*}(s^*)$ and $a_n=(f_{a^*,b^*}(s^*.y_n),n^*,\top)$ for some $y_n\geq x^*$ in $\mathbb{N}$, for any $k\in A_{\sqsubset_4}$ or $k\in A_{\sqsubset_{1};\sqsubset_{4}}$. As in  Case 1.2.1, we then obtain that $p_n=r_0$ is a constant, and so $a_n$ is fixed,  contradicting  the  assumption that $D_{\sqsubset_4}$ or $D_{\sqsubset_{1};\sqsubset_{4}}$ is infinite.

Case $1.2.4$: $D_{\sqsubset_1}$ is infinite. Then $a_n=(m^*,n^*,\top)$ for any $n\in A_{\sqsubset_1}$,
and so $a_n$ is fixed,  contradicting  the  assumption that $D_{\sqsubset_1}$ is infinite.

\textbf{Case $2$}: $(y_i)_{i\in I}\cap B$ is cofinal. Assume $(y_i)_{i\in I}\cap B=(m_0,n_0,(x_k)_{a_0,b_0})$ $_{k\in N}$ for some fixed natural numbers $m_0,n_0,a_0,b_0$, where $N$ is an infinite subset of $\mathbb{N}$. For each $k\in N$, there is $a_k\in \max B$ such that $(x_0,(m_0,n_0,(x_k)_{a_0,b_0}))< (a_k,a_k)$. Let $a_k=(m_k,n_k,\top)$ for every $k\in N$. By the order $<$ defined on $B$, we can see that for each $n$, either $n_k=n_0$ or $n_k=n_0+1$ holds.  Now consider two distinct cases for $x_0$.

\textbf{Case $2.1$}: $x_0\in P_1\backslash B$. Then assume $x_0=(f_0,r_0)$ for some fixed $(f_0,r_0)\in \mathbb{N}^{\mathbb{N}}\times \mathbb{N}$. Note that $x_0\leq a_k$ for all  $k\in N$. Then $x_0<_1;<_3 a_k$ or $x_0<_3 a_k$. It follows that $a_k=(f_0(z_k),z_k,\top)$ for some $z_k\geq r_0$, for each  $k\in N$. The infiniteness of $\{a_k:k\in N\}$ ensures that  $\{z_k:k\in N\}$ is infinite, a contradiction to that $z_k=n_k\in \{n_0,n_0+1\}$ for each $k\in N$ in Case 2.

\textcolor{black}{
\textbf{Case $2.2$}: $x_0\in B$. This is an symmetrical counterpart of the Case $2$ of Claim $2.1$, and we omit the argument. }

All the above cases show that Claim 2.2 is true. 

Finally, we have proved that $A$ is an irreducible closed subset of $\Sigma (P_1\times P_2)$. Note that $\sup A=(\top_1,\top_2)$, but $(\top_1,\top_2)\not\in A$. So $A$ is not a closure of some point in $P_1\times P_2$ and 
 hence,  $\Sigma (P_1\times P_2)$ cannot be  sober in the Scott topology. 
\end{proof}
\end{theorem}

\begin{remark}
Theorem~\ref{productisnotsober} shows that the category of all sober dcpo's and Scott-continuous maps is not closed under products in that of all dcpo's. Hence, 
the category of all sober dcpo's is not reflective in the category of all dcpo's and  Scott continuous mappings. 
\end{remark}

\section{Conclusion}
In this paper, we constructed two sober dcpo's, whose product is not sober. This gives a negative answer to a long-standing open problem in domain theory. As the two constructed dcpo's are not complete lattices, it is still not known whether the poset product of two sober complete lattices is sober.
We also proved that if $M=\sigma(P), N=\sigma(Q)$ with $P, Q$ being countable posets, then the poset product $M\times N$ is sober.

\bibliographystyle{plain}
\bibliography{refs}

@string{cup = {Cambridge University Press}}

@string{elsevier = {Elsevier Science Publishers {B.V.}}}

@string{ieee = {IEEE Computer Society Press}}

@string{is = {Information Systems}}

@string{lnm = {Lecture Notes in Mathematics}}

@string{sv = {Springer Verlag}}

@article{hochster1969prime,
  title={Prime ideal structure in commutative rings},
  author={Hochster, Melvin},
  journal={Transactions of the American Mathematical Society},
  volume={142},
  pages={43--60},
  year={1969}
}

@article{jia2020order,
  title={The order-sobrification monad},
  author={Jia, Xiaodong},
  journal={Applied Categorical Structures},
  volume={28},
  number={5},
  pages={845--852},
  year={2020},
  publisher={Springer}
}

@article{zhao2018uniqueness,
  title={Uniqueness of directed complete posets based on {Scott} closed set lattices},
  author={Zhao, Dongsheng and Xu, Luoshan},
  journal={Logical Methods in Computer Science},
  volume={14},
  year={2018},
  publisher={Episciences. org}
}

@article{xu2021some,
  title={Some open problems on well-filtered spaces and sober spaces},
  author={Xu, Xiaoquan and Zhao, Dongsheng},
  journal={Topology and its Applications},
  volume={301},
  pages={107540},
  year={2021},
  publisher={Elsevier}
}

@incollection{abramsky87a,
	Author = {Abramsky, Samson},
	Booktitle = {Symposium on Logic In Computer Science},
	Pages = {47--53},
	Publisher = ieee,
	Title = {Domain Theory in Logical Form},
	Year = {1987}}

@book{gierz03,
	Author = {Gierz, Gerhard and Hofmann, Karl H. and Keimel, Klaus and Lawson, Jimmie D. and Mislove, Michael and Scott, Dana S.},
	Publisher = cup,
	address = {Cambridge},
	Series = {Encyclopedia of Mathematics and its Applications},
	Title = {Continuous Lattices and Domains},
	Volume = 93,
	Year = {2003}}

@book{goubault13a,
	address = {Cambridge, UK},
	Author = {Goubault-Larrecq, Jean},
	Publisher = cup,
	Series = {New Mathematical Monographs},
	Title = {Non-Hausdorff Topology and Domain Theory},
	Volume = 22,
	Year = 2013}

@article{isbell82,
	Author = {John Isbell},
	Journal = {Proceedings of the American Mathematical Society},
	Pages = {333--334},
	Title = {Completion of a construction of {Johnstone}},
	Volume = 85,
	Year = 1982}

@inproceedings{johnstone81,
	Author = {Peter. T. Johnstone},
	Booktitle = {Continuous Lattices, Proceedings Bremen 1979},
	Editor = {B. Banaschewski and R.-E. Hoffmann},
	Pages = {282--283},
	Publisher = sv,
	Series = lnm,
	Title = {Scott is not always sober},
	Volume = {871},
	Year = {1981}}

@article{yukou15,
     Author = {Yue Yu and Hui Kou},
     Journal = {{Journal of Sichuan University (Natural Science Edition)}},
     Pages = {217-222},
     Title = {{D}irected spaces defined through {$T_0$} spaces with specialization order (in {Chinese})},
     Volume = {52(2)},
     Year = 2015}

@article{Miaoxl2023, 
title={Not every countable complete distributive lattice is sober}, 
volume={33}, 
DOI={10.1017/S0960129523000269}, 
number={9}, 
journal={Mathematical Structures in Computer Science}, 
author={Miao, Hualin and Xi, Xiaoyong and Li, Qingguo and Zhao, Dongsheng}, 
year={2023}, 
pages={809–831}}

@article{xuxizhao2021, 
title={A complete {Heyting} algebra whose {Scott} space is non-sober},
volume={252}, 
number={3}, 
journal={Fundamenta Mathematicae}, 
author={Xu, Xiaoquan and Xi, Xiaoyong and Zhao, Dongsheng}, 
year={2021}, 
pages={315–323}}

\end{document}